\documentclass[a4paper]{amsart}
\pdfoutput=1
\usepackage[T1]{fontenc}
\usepackage[utf8]{inputenc}
\usepackage{amssymb}
\usepackage{enumerate}
\usepackage{mathrsfs}

\usepackage[pdftitle={Combable groups have group cohomology of polynomial growth},
  pdfauthor={Ralf Meyer},
  pdfsubject={Mathematics; MSC 20F32, 20J05, 20J06}
]{hyperref}
\usepackage[lite]{amsrefs}
\usepackage{microtype}

\newcommand*{\brd}{-\hspace{0pt}}
\newcommand*{\nbd}{\nobreakdash-\hspace{0pt}}

\newcommand*{\abs}[1]{\lvert#1\rvert}          
\newcommand*{\norm}[1]{\lVert#1\rVert}         
\newcommand*{\trinorm}[1]{\lvert\!\lvert\!\lvert#1\rvert\!\rvert\!\rvert}
\newcommand*{\blank}{\text{\textvisiblespace}}
\newcommand*{\pt}{\star}                       

\newcommand*{\hot}{\mathbin{\hat{\otimes}}}    
\newcommand*{\defeq}{\mathrel{:=}}
\newcommand*{\din}{\precsim} 

\newcommand*{\Mod}{\mathsf{Mod}}

\newcommand*{\ID}{\mathrm{id}}      
\newcommand*{\pol}{\mathrm{pol}}    
\newcommand*{\Barr}{\mathrm{Bar}}   
\newcommand*{\Cell}{\mathrm{Cell}}  

\newcommand*{\T}{{\mathscr{T}}}     
\newcommand*{\Sch}{\mathscr{S}}     
\newcommand*{\Fil}{{\mathscr{F}}}   
\newcommand*{\Hol}{{\mathscr{O}}}   
\newcommand*{\Coh}{{\mathscr{C}}}   

\newcommand*{\R}{{\mathbb{R}}}
\newcommand*{\Z}{{\mathbb{Z}}}
\newcommand*{\N}{{\mathbb{N}}}
\newcommand*{\F}{{\mathbb{F}}}

\DeclareMathOperator{\Hom}{Hom}  
\DeclareMathOperator{\Ext}{Ext}  

\theoremstyle{plain}
\newtheorem{theorem}{Theorem}
\newtheorem{proposition}[theorem]{Proposition}
\newtheorem{lemma}[theorem]{Lemma}
\newtheorem{corollary}[theorem]{Corollary}
\theoremstyle{definition}
\newtheorem{definition}[theorem]{Definition}
\theoremstyle{remark}
\newtheorem{example}[theorem]{Example}

\begin{document}

\title{Combable groups have group cohomology of polynomial growth}
\author{Ralf Meyer}
\address{Mathematisches Institut\\
         Westfälische Wilhelms-Universität Münster\\
         Einsteinstr.\ 62\\
         48149 Münster\\
         Germany
}
\email{rameyer@math.uni-muenster.de}


\thanks{This research was supported by the EU-Network \emph{Quantum
    Spaces and Noncommutative Geometry} (Contract HPRN-CT-2002-00280)
  and the \emph{Deutsche Forschungsgemeinschaft} (SFB 478).}

\begin{abstract}
  Group cohomology of polynomial growth is defined for any finitely generated
  discrete group, using cochains that have polynomial growth with respect to
  the word length function.  We give a geometric condition that guarantees
  that it agrees with the usual group cohomology and verify this condition for
  a class of combable groups.  Our condition involves a chain complex that is
  closely related to exotic cohomology theories studied by Allcock and Gersten
  and by Mineyev.
\end{abstract}
\maketitle

\section{Introduction}
\label{sec:intro}

Let~\(G\) be a finitely generated discrete group and let~\(\ell\) be a word
length function on~\(G\).  In a standard complex that computes the group
cohomology of~\(G\), one can define a subcomplex by allowing only cochains of
polynomial growth with respect to~\(\ell\) in each variable.  The cohomology
of this subcomplex is the \emph{group cohomology of polynomial growth}
of~\(G\), which we denote by \(H^n_\pol(G)\).  We may ask whether the
canonical maps \(H^n_\pol(G)\to H^n(G)\) are isomorphisms.  This question came
up in the work of Alain Connes and Henri Moscovici on the Novikov conjecture
for hyperbolic groups in~\cite{Connes-Moscovici:Novikov_Hyperbolic} and has
also been studied by Ronghui Ji in~\cite{Ji:Schwartz_cohomology}.

We define a certain chain complex of bornological vector spaces
\(\Sch\tilde{C}_\bullet(G)\) and prove the following two facts.  First,
\(\Sch\tilde{C}_\bullet(G)\) has a bounded contracting homotopy if the
group~\(G\) is combable in the sense of~\cite{Epstein:Automatic}.  Secondly,
the existence of such a contracting homotopy implies \(H^n_\pol(G)\cong
H^n(G)\).  An important feature of our construction is that the homotopy type
of the chain complex \(\Sch\tilde{C}_\bullet(G)\) is a quasi\brd{}isometry
invariant of~\(G\).  In contrast, the problem of whether \(H^n_\pol(G)\cong
H^n(G)\) does not seem invariant under quasi\brd{}isometry.

The chain complex \(\Sch\tilde{C}_\bullet(G)\) is related to the convolution
algebra
\[
\Sch(G) \defeq
\Bigl\{ f\colon G\to\R \Bigm|
  \sum_{g\in G} \abs{f(g)} (\ell(g)+1)^k<\infty \quad \forall k\in\N
\Bigr\}.
\]
If we use the Banach algebra \(\ell_1(G)\) instead, we obtain a complex that
is homotopy equivalent to one constructed by Daniel J.\ Allcock and Stephen
M.\ Gersten in~\cite{Allcock-Gersten} (for groups for which the classifying
space \(BG\) has finite type).  However, already the chain complex
\(\ell_1\tilde{C}_\bullet(\Z)\) has non-trivial homology and hence cannot be
contractible.  The best one can say is that the range of the differential in
this complex is dense in the kernel, that is, the ``reduced homology''
vanishes.  Igor Mineyev shows in~\cite{Mineyev:Homology_Combable} that this
happens for groups with a sufficiently nice combing.  However, this seems too
weak to compare the bounded cohomology and the usual group cohomology.

We now explain the definition of \(\Sch\tilde{C}_\bullet(G)\).  For any
set~\(X\), let \(S(X)\) be the simplicial set whose \(n\)\nbd{}simplices are
the \({n+1}\)\brd{}tuples \((x_0,\dotsc,x_n)\in X^{n+1}\) and whose \(j\)th
face and degeneracy maps leave out and double~\(x_j\), respectively.  Let
\(C_\bullet(X)\) be the reduced simplicial chain complex associated to
\(S(X)\) with coefficients~\(\R\).  We identify \(C_n(X)\) with the space of
functions \(X^{n+1}\to\R\) with finite support that vanish on
\((x_0,\dotsc,x_n)\) if \(x_j=x_{j+1}\) for some \(j\in\{0,\dotsc,n-1\}\).  We
define the augmentation map \(\alpha\colon C_0(X)\to\R\) on basis vectors by
\(x\mapsto 1\) for all \(x\in X\).  We define a subcomplex
\(\tilde{C}_\bullet(X)\) by \(\tilde{C}_0(X) = \ker \alpha\) and
\(\tilde{C}_n(X)=C_n(X)\) for \(n\ge1\).

Any map \(f\colon X\to Y\) induces a chain map \(f_*\colon
\tilde{C}_\bullet(X)\to\tilde{C}_\bullet(Y)\), and there is an explicit chain
homotopy between \(f_*\) and~\(f'_*\) for any \(f,f'\colon X\to Y\).  Hence
\(\tilde{C}_\bullet(X)\) is always contractible because it is homotopy
equivalent to \(\tilde{C}_\bullet(\pt) = 0\).

If \((X,d)\) is a discrete proper metric space, we let \(\Sch C_n(X)\) be the
space of functions \(f\colon X^{n+1}\to\R\) with the following properties:
\begin{enumerate}[(1)]
\item there is \(R>0\) such that \(f(x_0,\dotsc,x_n)=0\) if \(d(x_i,x_j)\ge
  R\) for some \(i,j\in\{0,\dotsc,n\}\), or if \(x_j=x_{j+1}\) for some
  \(j\in\{0,\dotsc,n-1\}\);
\item for all \(k\in\N\) and all fixed \(\pt\in X\), we have
  \[
  \sum_{x_0,\dotsc,x_n\in X} \abs{f(x_0,\dotsc,x_n)} \cdot
  (d(x_0,\pt)+\dotsb+d(x_n,\pt)+1)^k
  < \infty.
  \]
\end{enumerate}
The space \(\Sch C_n(X)\) is a bornological vector space that contains
\(C_n(X)\) as a dense subspace.  The differential of \(C_\bullet(X)\) may be
extended and turns \(\Sch C_\bullet(X)\) into a chain complex of bornological
vector spaces.  Finally, \(\Sch\tilde{C}_\bullet(X)\) is the kernel of the
augmentation map \(\Sch C_\bullet(X)\to\R\).

A combing on a metric space \((X,d)\) with a chosen base point~\(\pt\) is a
sequence of maps \(f_n\colon X\to X\) with the following properties:
\(f_0(x)=\pt\) for all \(x\in X\), and for any \(x\in X\) there exists
\(n\in\N\) such that \(f_N(x)=x\) for all \(N\ge n\); the maps~\(f_n\) are
uniformly quasi\brd{}Lipschitz; and the pairs of maps \((f_n,f_{n+1})\) are
uniformly close in the sense that the set of \(d(f_n(x),f_{n+1}(x))\) for
\(x\in X\), \(n\in\N\) is bounded.  We say that the combing has polynomial
growth if the number of \(n\in\N\) with \(f_n(x)\neq f_{n+1}(x)\) is
controlled by a polynomial in \(d(x,\pt)\).  Such combings exist for all
hyperbolic groups, for automatic groups, and for groups that are combable in
the sense of~\cite{Epstein:Automatic}.  However, it seems that such combings
do not exist for non-Abelian nilpotent groups.

Our first main result is that \(\Sch\tilde{C}_\bullet(X)\) has a bounded
contracting homotopy if~\(X\) has a combing of polynomial growth.  The proof
is actually quite simple.  Since the maps \(f_n,f_{n+1}\) are close for each
\(n\in\N\), there is an explicit chain homotopy \(H(f_n,f_{n+1})\colon
\Sch\tilde{C}_\bullet(X)\to \Sch\tilde{C}_\bullet(X)\) between the maps
induced by \(f_n\) and~\(f_{n+1}\).  The sum \(H = \sum_{n\in\N}
H(f_n,f_{n+1})\) is the desired contracting homotopy.  The hypotheses on
\((f_n)\) guarantee that~\(H\) is a bounded linear operator on
\(\Sch\tilde{C}_\bullet(X)\).

Our second main result is that \(H^n_\pol(G)\cong H^n(G)\) if
\(\Sch\tilde{C}_\bullet(G)\) has a bounded contracting homotopy.  This
requires some homological algebra with bornological modules over the
convolution algebras \(\R[G]\) and \(\Sch(G)\).  These two are bornological
unital algebras, and the embedding \(i\colon \R[G]\to\Sch(G)\) is a bounded
unital algebra homomorphism.  We work with bornologies instead of topologies
because this gives better results for spaces that are built out of \(\R[G]\)
and \(\Sch(G)\).  The right category of modules over a bornological unital
algebra~\(A\) (like \(\R[G]\) or \(\Sch(G)\)) is the category \(\Mod(A)\) of
bornological left \(A\)\nbd{}modules.  Such a module is defined by a bounded
unital algebra homomorphism \(A\to\Hom(M,M)\) or, equivalently, a bounded
bilinear map \(A\times M\to M\) satisfying the usual properties.  The
morphisms are the bounded \(A\)\nbd{}linear maps.  Homological algebra in
\(\Mod(A)\) works as usual if we only admit extensions with a bounded linear
section and resolutions with a bounded linear contracting homotopy.  This is
explained in detail in~\cite{Meyer:Embed_derived}.

In order to keep this article accessible for readers with a background in
geometric group theory, we only use the most basic homological algebra and do
not explore all consequences of the contractibility of
\(\Sch\tilde{C}_\bullet(G)\).  This is done in~\cite{Meyer:Embed_derived},
where I also show by other techniques that \(\Sch\tilde{C}_\bullet(G)\) is
contractible for groups of polynomial growth.  This improves the main result
of~\cite{Ji:Schwartz_cohomology}.

Let~\(\hot\) be the projective complete bornological tensor product
(\cite{Hogbe-Nlend:Completions}).  A bornological \(A\)\nbd{}module of the
form \(A\hot X\) with the obvious module structure is called \emph{free}.
Free modules are projective for extensions with a bounded linear section.  The
usual argument in homological algebra that shows that two free resolutions of
the same module are homotopy equivalent still works in our setting.

The group cohomology of~\(G\) can be defined as \(H^n(G) \cong
\Ext_{\R[G]}^n(\R,\R)\), where~\(\R\) is equipped with the trivial
representation of~\(G\) and the resulting module structure over \(\R[G]\).
That is, it is the cohomology of the chain complex
\(\Hom_{\R[G]}(P_\bullet,\R)\), where \(P_\bullet\to\R\) is some free
\(\R[G]\)\brd{}module resolution of~\(\R\).  If we let~\(G\) act diagonally on
\(C_n(G)\), then the chain complex \(C_\bullet(G)\) constructed above becomes
such a free \(\R[G]\)\brd{}module resolution.  Similarly, the group cohomology
of polynomial growth is isomorphic to the cohomology
\(\Ext^n_{\Sch(G)}(\R,\R)\) of the chain complex
\(\Hom_{\Sch(G)}(P'_\bullet,\R)\), where \(P'_\bullet\to\R\) is some free
\(\Sch(G)\)\brd{}module resolution of~\(\R\).

Let \(i^*\colon \Mod(\Sch(G))\to\Mod(\R[G])\) be the functor induced by the
embedding \(i\colon \R[G]\to\Sch(G)\).  We also define a functor
\[
i_!\colon \Mod(\Sch(G))\to\Mod(\R[G]),
\qquad i_!(M) \cong \Sch(G)\hot_{\R[G]} M.
\]
These functors are adjoint, that is,
\begin{equation}
  \label{eq:istar_adjoint}
  \Hom_{\Sch(G)}(i_!(M),N) \cong \Hom_{\R[G]}(M,i^*(N))
\end{equation}
if \(M\) and~\(N\) are bornological left modules over \(\R[G]\) and
\(\Sch(G)\), respectively.  Moreover, \(i_!(C_\bullet(G)) \cong \Sch
C_\bullet(G)\).  If \(\Sch\tilde{C}_\bullet(G)\) has a bounded contracting
homotopy, then \(i_!(C_\bullet(G))\) is a free \(\Sch(G)\)\brd{}module
resolution of~\(\R\).  Hence we may use it to compute
\(\Ext^n_{\Sch(G)}(\R,\R)\).  Since \(i^*(\R)=\R\), the adjointness
relation~\eqref{eq:istar_adjoint} yields
\begin{multline*}
  \Ext^n_{\Sch(G)}(\R,\R)
  \cong H^n(\Hom_{\Sch(G)}(i_!(C_\bullet(G)),\R))
  \\ \cong H^n(\Hom_{\R[G]}(C_\bullet(G),\R))
  \cong \Ext^n_{\R[G]}(\R,\R).
\end{multline*}

We now summarise the following sections.  In Section~\ref{sec:function_spaces}
we define bornological vector spaces \(\Sch(X)=\Sch^\infty(X)\),
\(\Sch^\omega(X)\) and \(\Hol(X)\) of functions of polynomial, subexponential,
and exponential decay.  Along the way, we recall some basic facts about
bornologies because they went out of fashion a few decades ago.  Let \(\T(X)\)
be one of the function spaces mentioned above.  In
Section~\ref{sec:def_chains_metric} we define a chain complex \(C_\bullet(X)\)
for any set~\(X\) and enlarge it to \(\T C_\bullet(X)\) if~\(X\) is a discrete
proper metric space.  We check that \(\T C_\bullet(X)\) is functorial for
quasi\brd{}Lipschitz maps and that close maps induce chain homotopic chain
maps.  In Section~\ref{sec:combings_contract} we explain the notion of combing
that we use and prove that \(\T C_\bullet(X)\) is contractible if~\(X\) has a
combing of appropriate growth.

In Section~\ref{sec:tempered_cohomology}, we define group cohomology with
polynomial, subexponential, and exponential growth and identify these theories
with extension groups in the bornological module categories of
\(\Sch^\infty(G)\), \(\Sch^\omega(G)\), and \(\Hol(G)\), respectively.  Then
we prove that the tempered group cohomology associated to \(\T(G)\) is
isomorphic to the usual group cohomology if the reduced chain complex
\(\T\tilde{C}_\bullet(G)\) has a bounded contracting homotopy.  In addition,
we prove that \(H^n(G)\) is finite-dimensional for all \(n\in\N\) if
\(\T\tilde{C}_\bullet(G)\) is contractible.  Hence
\(\Sch\tilde{C}_\bullet(G)\) may fail to be contractible.

It is desirable to replace \(\T C_\bullet(X)\) by smaller and more geometric
chain complexes.  We discuss how to do this in Sections
\ref{sec:other_resolution} and~\ref{sec:tempered_cellular_chains}.  If~\(X\)
is a group, then \(\T C_\bullet(X)\) is homotopy equivalent to
\(i_!(P_\bullet)\) for any free \(\R[G]\)\brd{}module resolution of~\(\R\).
This allows to study the case of non-Abelian free groups in great detail.  One
particular way to construct free \(\R[G]\)\brd{}module resolutions of the
trivial representation is to use the cellular chain complex
\(\Cell_\bullet(Y)\) of a contractible CW\brd{}complex on which~\(G\) acts
properly by cellular maps.  We investigate what additional structure on~\(Y\)
is needed to write down the tempered versions of \(\Cell_\bullet(Y)\).  This
yields a sufficient condition for \(\T\Cell_\bullet(Y)\cong \T C_\bullet(X)\)
that works if~\(X\) is just a proper discrete metric space.

The results of this article were announced at a conference at the
Mathematische Forschungsinstitut Oberwolfach in 2004 (see
\cite{Meyer:Poly_Oberwolfach}).  More recently, I have applied similar
techniques to the Schwartz algebras of reductive \(p\)\nbd{}adic groups
(\cite{Meyer:Homological_Schwartz}).

I would like to thank Heath Emerson for helpful discussions about coarse
geometry and geometric group theory.

\section{Function spaces on discrete metric spaces}
\label{sec:function_spaces}

Let~\(X\) be a discrete, proper metric space.  Here and in the following,
\emph{properness} means that bounded subsets of~\(X\) are finite.  Let
\(\R[X]\) be the free \(\R\)\nbd{}vector space on the set~\(X\); that is,
\(\R[X]\) is a vector space containing~\(X\) as a basis.  We often identify
\(\R[X]\) with the space of functions \(X\to\R\) of finite support.  Fix any
point \(\pt\in X\) and define \(\ell(x)\defeq d(x,\pt)\).  Define norms on
\(\R[X]\) by
\[
\norm{\phi}^k \defeq \sum_{x\in X} \abs{\phi(x)} (\ell(x)+1)^k
\]
for any \(k\in\R\).  We let \(\Sch^k(X)\) be the Banach space completion of
\(\R[X]\) for this norm.  Equivalently, \(\Sch^k(X)\) consists of those
functions \(\phi\colon X\to\R\) with \(\norm{\phi}^k<\infty\).  We get
equivalent norms if we change~\(\pt\) or replace~\(d\) by a metric that is
quasi\brd{}isometric to it (that is, the identity map \((X,d)\to (X,d')\) is a
quasi\brd{}isometry).  Hence such modifications yield equivalent spaces
\(\Sch^k(X)\).  Let \(\Sch(X)\) be the space of functions \(\phi\colon
X\to\R\) that satisfy \(\norm{\phi}^k<\infty\) for all \(k\in\N\).  This is a
Fréchet space with respect to the topology that is defined by the sequence of
norms \(\norm{\blank}^k\), \(k\in\N\).  Later we shall usually treat
\(\Sch(X)\) as a bornological vector space.

\smallbreak

A \emph{bornology} on a vector space is a collection of subsets, called
\emph{bounded} subsets, satisfying the following axioms
(see~\cite{Hogbe-Nlend:Bornologies}): subsets of bounded subsets are again
bounded; finite subsets are bounded; if \(S_1,S_2\) are bounded, so is
\(S_1+S_2\); if~\(S\) is bounded, so is \(\bigcup_{\abs{t}\le\lambda} t\cdot
S\) for any \(\lambda\in\R\).
Two important examples are the families of von Neumann bounded subsets and of
precompact subsets of a Fréchet space.  We call these bornologies the
\emph{von Neumann bornology} and the \emph{precompact bornology}.  Recall that
a subset of a Fréchet space~\(V\) is called \emph{von Neumann bounded} if it
is bounded with respect to all continuous seminorms.  Since~\(V\) is complete,
there is no difference between relatively compact and precompact subsets.  By
the way, ``precompact'' is a synonym for ``totally bounded''.  Usually, the
von Neumann bornology is easier to describe explicitly than the precompact
bornology, but the precompact bornology may have better analytical properties.

One checks easily that the embedding \(\Sch^{k+1}(X)\to\Sch^k(X)\) is a
compact map for all \(k\in\N\).  By definition, this means that \(\Sch(X)\) is
a \emph{Fréchet-Schwartz space}.  In a Fréchet-Schwartz space, any von Neumann
bounded subset is relatively compact, that is, the von Neumann bornology and
the precompact bornology agree (see \cite{Hogbe-Nlend:Bornologies}).  We
always equip \(\Sch(X)\) with this bornology.  We may equip the Banach spaces
\(\Sch^k(X)\) with the precompact bornology or the von Neumann bornology.
This makes no difference for the following constructions.

In this article all bornologies are complete and convex.  This means that for
any bounded subset \(S\subseteq V\) there exists a subspace \(V_T\subseteq V\)
together with a Banach space norm on~\(V_T\) such that the embedding \(V_T\to
V\) is a bounded linear map and such that~\(S\) is a bounded subset
in~\(V_T\).  The Banach space norm on~\(V_T\) can be described equally well by
its closed unit ball.  This is a subset \(T\subseteq V\), which we call a
\emph{complete disk} in~\(V\).  The complete disks in any complete
bornological vector space~\(V\) form a directed set.  The space~\(V\) is
canonically isomorphic to the direct limit of the associated inductive system
of Banach spaces \((V_T)\).  Various concepts of functional analysis are
defined in the bornological context by reduction to these Banach spaces.  For
instance, a sequence in~\(V\) converges if and only if it converges in~\(V_T\)
for some complete disk~\(T\).  Since we only use complete convex bornologies
in the following, we drop these qualifiers from our notation and agree that
bornology, henceforth, stands for complete convex bornology.

Another important example of a bornology is the \emph{fine bornology}, which
is defined for any vector space~\(V\).  It consists of those subsets that are
contained in and bounded in some finite-dimensional subspace of~\(V\).  It is
the finest (that is, smallest) possible bornology on~\(V\).  This is a
reasonable bornology if~\(V\) is a vector space like \(\R[X]\).  Analysis in
fine bornological vector spaces is quite easy because it always reduces to
finite-dimensional subspaces.  In contrast, the finest locally convex topology
on a space like \(\R[X]\) is rather unwieldy.  The bornological and
topological ways of doing analysis turn out to be equivalent for Fréchet
spaces like \(\Sch(X)\) in many situations (see~\cite{Meyer:Born_Top}).  For
spaces like \(\R[X]\) there are some technical differences, where the
bornological approach gives better results.

A linear map \(f\colon V\to W\) between two bornological vector spaces is
called \emph{bounded} if it maps bounded subsets of~\(V\) to bounded subsets
of~\(W\).  The bounded linear maps are the morphisms in the category of
bornological vector spaces.  If \(V\) and~\(W\) are both Fréchet spaces
equipped with the von Neumann or precompact bornology, then a linear map
\(V\to W\) is bounded if and only if it is continuous.  If~\(V\) carries the
fine bornology and~\(W\) is arbitrary, then any linear map \(V\to W\) is
bounded.

We frequently have to express the fact that some bornological vector
space~\(V\), say, \(V=\R[X]\) is contained in another bornological vector
space~\(W\), say, \(W=\Sch(X)\), such that the embedding \(V\to W\) is a
bounded linear map with dense range.  We simply write \(V\din W\) to express
all of the above.  ``Dense range'' means that a bounded linear map \(W\to X\)
into another bornological vector space that vanishes on~\(V\) already vanishes
on all of~\(W\) (see also~\cite{Meyer:Born_Top} for other notions of density).

\smallbreak

After this digression into functional analysis, we construct some more
function spaces.  For any \(\alpha>0\), let \(\ell_1(X,\alpha^\ell) \defeq
\{\phi\colon X\to\R\mid \trinorm{\phi}_\alpha<\infty\}\), where
\[
\trinorm{\phi}_\alpha \defeq \sum_{x\in X} \abs{\phi(x)} \alpha^{\ell(x)}.
\]
Although \(\trinorm{\blank}_\alpha\) and \(\ell_1(X,\alpha^\ell)\) change if we
rescale the metric, the spaces
\begin{align*}
  \Hol(X) &\defeq \{\phi\colon X\to\R\mid
    \forall \alpha>1 \colon\quad \trinorm{\phi}_\alpha < \infty
  \}
  \\
  \Sch^\omega(X) &\defeq \{\phi\colon X\to\R\mid
    \exists \alpha>1 \colon\quad \trinorm{\phi}_\alpha < \infty
  \}
\end{align*}
only depend on the quasi\brd{}isometry class of the metric.

The sequence of norms \(\trinorm{\blank}_n\) for \(n\in\N\) turns \(\Hol(X)\)
into a Fréchet-Schwartz space because the embeddings
\(\iota_{\alpha,\beta}\colon \ell_1(X,\alpha^\ell)\to \ell_1(X,\beta^\ell)\)
for \(\alpha>\beta\) are compact.  Hence the von Neumann bornology and the
precompact bornology on \(\Hol(X)\) agree.  We choose this bornology on
\(\Hol(X)\).  We equip \(\Sch^\omega(X)\) with the direct limit bornology.
Thus a subset of \(\Sch^\omega(X)\) is bounded if and only if it is bounded
with respect to the norm \(\trinorm{\blank}_{1+1/n}\) for some \(n\in\N\).  A
bornological vector space that is a direct limit of a sequence of Banach
spaces with compact structure maps is called a \emph{Silva space}
(see~\cite{Hogbe-Nlend:Bornologies}).  Thus \(\Sch^\omega(X)\) is a Silva
space.

There is a duality between Fréchet-Schwartz spaces and Silva spaces,
see~\cite{Hogbe-Nlend:Bornologies}.  The dual of a topological vector space
carries a canonical bornology: the equicontinuous bornology, consisting of
equicontinuous families of linear functionals.  The dual of a Fréchet-Schwartz
space with this bornology is a Silva space.  The dual of a bornological vector
space carries a canonical topology: the topology of uniform convergence on
bounded subsets.  The dual of a Silva space with this topology is a
Fréchet-Schwartz space.  These two operations are inverse to each other.

We need to know the dual spaces of the function spaces defined above.  Since
\(\R[X]\) is a dense subspace in any of them, a bounded linear functional is
determined by its values on the basis vectors \(x\in X\) and hence given by
\(\phi\mapsto \sum_{x\in X} \phi(x) \psi(x)\) for some function \(\psi\colon
X\to\R\).  Fix such a function~\(\psi\).  It gives rise to a bounded linear
functional on \(\Sch^k(X)\) for \(k\in\N\) if and only if~\(\psi\) has
\emph{polynomial growth of order~\(k\)}, that is, \(\psi(x)\le
C(\ell(x)+1)^k\) for some \(C>0\), \(k\in\N\).  It gives rise to a bounded
linear functional on \(\Sch(X)\) if and only if it has \emph{polynomial
  growth} of some order \(k\in\N\).  It gives rise to a bounded linear
functional on \(\Sch^\omega(X)\) if and only if it has \emph{subexponential
  growth}, that is, for any \(\alpha>1\) there is \(C>0\) such that
\(\psi(x)\le C\cdot \alpha^{\ell(x)}\) for all \(x\in X\).  A
function~\(\psi\) gives rise to a bounded linear functional on \(\Hol(X)\) if
and only if it has \emph{exponential growth}, that is, \(\psi(x)\le C\cdot
\alpha^{\ell(x)}\) for all \(x\in X\) for some \(C>0\), \(\alpha>1\).  All
these dual spaces carry canonical bornologies.  The spaces of functions of
polynomial and exponential growth are Silva spaces, the space of functions of
subexponential growth is a Fréchet-Schwartz space.

Let \(X\) and~\(Y\) be discrete proper metric spaces.  A map \(f\colon X\to
Y\) is called \emph{quasi\brd{}Lipschitz} or \emph{large scale Lipschitz} if
there is \(C>0\) such that \(d(f(x),f(y))\le C(d(x,y)+1)\) for all \(x,y\in
X\).  Then the linear map \(f_*\colon \R[X]\to\R[Y]\) defined by
\(f_*(\delta_x)\defeq \delta_{f(x)}\) or, equivalently, \(f_*\phi(y) \defeq
\sum_{f(x)=y} \phi(x)\), extends to bounded linear operators
\(\Sch^k(X)\to\Sch^k(Y)\) for \(k\in\R_+\cup\{\infty,\omega\}\) and
\(\Hol(X)\to\Hol(Y)\).  This is a special feature of \(\ell_1\)\nbd{}norms.
For instance, \(\ell_2(X)\) is only functorial if the number of pre-images of
points in~\(Y\) is bounded.  Thus \(\ell_2(X)\) is not even functorial for
quasi\brd{}isometries.

Now let~\(G\) be a finitely generated discrete group equipped with a word
metric, so that~\(\ell\) is the word length function associated to some set of
generators.  Hence \(\ell(xy)\le\ell(x)+\ell(y)\) for all \(x,y\in G\).  This
implies easily that the norms \(\norm{\blank}^k\) and
\(\trinorm{\blank}_\alpha\) are submultiplicative.  Thus \(\Sch^k(G)\) and
\(\ell_1(G,\alpha^\ell)\) are Banach algebras for all \(k\in\N\),
\(\alpha>1\).  It follows that \(\Sch^\infty(G)\) and \(\Hol(G)\) are locally
multiplicative Fréchet algebras and that \(\Sch^\omega(G)\) is a bornological
algebra.

\section{Tempered free chain complexes}
\label{sec:def_chains_metric}

We first define a chain complex \(C_\bullet(X)\) for any non-empty set~\(X\)
and discuss its functoriality.  Then we modify it using a metric on~\(X\).

The \emph{free simplicial set} \(S(X)\) over a non-empty set~\(X\) is the
following simplicial set.  We let \(S_n(X)\) be the set of all
\({n+1}\)\brd{}tuples \((x_0,\dotsc,x_n)\in X^{n+1}\).  The \(j\)th face map
\(S_n(X)\to S_{n-1}(X)\) omits~\(x_j\), the \(j\)th degeneracy map \(S_n(X)\to
S_{n+1}(X)\) doubles~\(x_j\).  Let \(C_n(X)\) be the vector space spanned by
the nondegenerate simplices in \(S_n(X)\).  Thus a basis of \(C_n(X)\) is
given by the \({n+1}\)\brd{}tuples \((x_0,\dotsc,x_n)\) with \(x_j\neq
x_{j+1}\) for all \(j\in\{0,\dotsc,n-1\}\).  Equivalently, \(C_n(X)\) is the
space of functions \(\phi\colon X^{n+1}\to\R\) with finite support that
satisfy \(\phi(x_0,\dotsc,x_n)=0\) if \(x_j=x_{j+1}\) for some
\(j\in\{0,\dotsc,n-1\}\).  We equip \(C_n(X)\) with the fine bornology.

The simplicial boundary map \(\delta\colon C_n(X)\to C_{n-1}(X)\) is given on
basis vectors by
\[
\delta(x_0,\dotsc,x_n)
\defeq \sum_{j=0}^n (-1)^j (x_0,\dotsc,\widehat{x_j},\dotsc,x_n).
\]
As usual, \(\widehat{x_j}\) means omission of~\(x_j\).  It is easy to check by
hand that~\(\delta\) is well-defined and satisfies \(\delta^2=0\).  In terms
of functions, we have
\[
\delta \phi(x_0,\dotsc,x_{n-1}) = \sum_{j=0}^n \sum_{y\in X} (-1)^j
\phi(x_0,\dotsc,x_{j-1},y,x_j,\dotsc,x_{n-1})
\]
for all \(\phi\in C_n(X)\).  We call \(C_\bullet(X) \defeq
(C_n(X),\delta)_{n\ge0}\) the \emph{free chain complex over~\(X\)}.

Any map \(f\colon X\to Y\) induces a morphism of simplicial sets \(S(X)\to
S(Y)\) by
\[
f_*(x_0,\dotsc,x_n)\defeq \bigl(f(x_0),\dotsc,f(x_n)\bigr)
\qquad \forall n\in\N,\ x_0,\dotsc,x_n\in X.
\]
This induces a chain map \(f_*\colon C_\bullet(X)\to C_\bullet(Y)\) by the
same formula.

Let \(f,f'\colon X\to Y\) be two maps.  There is an elementary homotopy
between the induced maps \(f_*,f'_*\colon S(X)\to S(Y)\).  On the chain
complex level, it is given by
\begin{multline*}
  H(f,f') \colon C_n(X)\to C_{n+1}(Y),
  \\
  (x_0,\dotsc,x_n) \mapsto \sum_{j=0}^n (-1)^j
  \bigl(f(x_0),\dotsc,f(x_j),f'(x_j),\dotsc,f'(x_n)\bigr).
\end{multline*}
Since~\(x_j\) occurs twice,
\(\bigl(f(x_0),\dotsc,f(x_j),f'(x_j),\dotsc,f'(x_n)\bigr)\) contains a
repetition if \((x_0,\dotsc,x_n)\) does.  Hence \(H(f,f')\) is well-defined.
We claim that
\begin{equation}
  \label{eq:Hff_chain_homotopy}
  H(f,f') \circ \delta + \delta \circ H(f,f') = f'_* - f_*.
\end{equation}
We verify~\eqref{eq:Hff_chain_homotopy} on a generator \((x_0,\dotsc,x_n)\).
All the summands in \(H(f,f') \delta\) also occur in \(\delta H(f,f')\) with
opposite signs and hence cancel out.  The only terms in \((-1)^j
\delta\bigl(f(x_0),\dotsc,f(x_j),f'(x_j),\dotsc,f'(x_n)\bigr)\) that survive
are
\[
\bigl(f(x_0),\dotsc,f(x_{j-1}),f'(x_j),\dotsc,f'(x_n)\bigr)
-\bigl(f(x_0),\dotsc,f(x_j),f'(x_{j+1}),\dotsc,f'(x_n)\bigr).
\]
Summing over~\(j\), another cancellation occurs and
yields~\eqref{eq:Hff_chain_homotopy}.

As a consequence, the constant maps \(X\to\pt\to X\) yield a homotopy
equivalence between \(S(X)\) and \(S(\pt)\) and hence between \(C_\bullet(X)\)
and \(C_\bullet(\pt)\) for any point \(\pt\in X\).  Evidently,
\(\abs{S(\pt)}=\pt\) and \(C_\bullet(\pt)\) is just~\(\R\) in degree~\(0\).
The homotopy equivalence \(C_\bullet(X)\to\R\) is the augmentation map
\begin{equation}
  \label{eq:augmentation_map}
  \alpha\colon C_0(X)= \R[X]\to\R,
  \qquad \alpha \phi = \sum_{x\in X} \phi(x).
\end{equation}
We let \(\tilde{C}_\bullet(X) \subseteq C_\bullet(X)\) be the subcomplex with
\(\tilde{C}_n(X)\defeq C_n(X)\) for \(n\ge1\) and \(\tilde{C}_0(X)\defeq
\ker\alpha \subseteq C_0(X)\).  Thus \(\tilde{C}_\bullet(X)\) is contractible
for any set~\(X\).

\smallskip

Now let~\(X\) be a discrete proper metric space.  Let \(\T(X)\) be one of the
function spaces \(\Sch^k(X)\) for \(k\in\R_+\cup\{\infty,\omega\}\) or
\(\Hol(X)\) that we have defined in Section \ref{sec:function_spaces}.  For
any \(R\in\N\), we let \(\T C_n(X)_R\) be the space of all functions
\(\phi\colon X^{n+1}\to\R\) that satisfy the following two conditions:
\begin{description}
\item[control] \(\phi(x_0,\dotsc,x_n)=0\) if \(d(x_i,x_j)>R\) for some
  \(i,j\in\{0,\dotsc,n\}\) or \(x_j=x_{j+1}\) for some
  \(j\in\{0,\dotsc,n-1\}\);
    
\item[growth] the function \(x\mapsto \sum_{x_1,\dotsc,x_n\in X}
  \abs{\phi(x,x_1,\dotsc,x_n)}\) belongs to \(\T(X)\).

\end{description}
A subset of \(\T C_n(X)_R\) is bounded if the set of functions in \(\T(X)\)
that arises in the second condition is bounded in \(\T(X)\).

Since bounded subsets of~\(X\) are finite, the control condition for functions
in \(\T C_n(X)_R\) ensures that the sums in the growth condition are finite.
The norms that define \(\T(X)\) only depend on~\(\ell\), and \(\abs{\ell(x) -
  \ell(y)}\le R\) if \(d(x,y)\le R\).  Hence we can rewrite the growth
condition more symmetrically.  A function \(\phi\colon X^{n+1}\to\R\) with
controlled support belongs to \(\Sch^k C_n(X)_R\) for \(k\in\R_+\) if and only
if
\[
\sum_{x_0,\dotsc,x_n\in X} \abs{\phi(x_0,\dotsc,x_n)}
(\ell(x_0)+\dotsb+\ell(x_n)+1)^k <\infty.
\]
We have \(\phi\in\Sch^\infty C_n(X)_R\) if and only if these norms remain
finite for all \(k\in\N\).  We leave it to the reader to formulate similar
descriptions for \(\Sch^\omega(X)_R\) and \(\Hol(X)_R\).

We let \(\T C_n(X)\) be the inductive limit of the subspaces \(\T C_n(X)_R\),
\(R\in\N\), equipped with the inductive limit bornology.  That is, a subset of
\(\T C_n(X)\) is bounded if it is contained in and bounded in \(\T C_n(X)_R\)
for some \(R\in\N\).  For \(n=0\) the control condition is empty and we simply
get \(\T C_0(X)= \T(X)\).  The augmentation map~\eqref{eq:augmentation_map}
extends to a bounded linear functional on \(\T(X)\).  Hence we get a
\(1\)\nbd{}codimensional closed subspace \(\T\tilde{C}_0(X)\defeq \ker
\alpha\subseteq \T C_0(X)\).  We also let \(\T\tilde{C}_n(X)=\T C_n(X)\) for
\(n\ge1\).  One checks easily that the boundary map~\(\delta\) extends to a
bounded linear operator \(\T\tilde{C}_n(X)\to \T\tilde{C}_{n-1}(X)\) for all
\(n\ge1\).  Furthermore, we have \(\tilde{C}_n(X)\din\T\tilde{C}_n(X)\) and
\(C_n(X)\din\T C_n(X)\) for all \(n\in\N\) because \(\R[X]\din \T(X)\).

\begin{definition}
  \label{def:tempered_free_chains}
  Let~\(X\) be a discrete proper metric space and let \(\T(X)\) be
  \(\Sch^k(X)\) for \(k\in\R_+\cup\{\infty,\omega\}\) or \(\Hol(X)\).  The
  chain complexes of bornological vector spaces
  \[
  \T C_\bullet(X) = (\T C_n(X),\delta)_{n\ge0},
  \qquad
  \T\tilde{C}_\bullet(X) = (\T\tilde{C}_n(X),\delta)_{n\ge0}
  \]
  are called \emph{\(\T(X)\)\brd{}tempered free chain complex over~\(X\)} and
  \emph{\(\T(X)\)\brd{}tempered reduced free chain complex over~\(X\)},
  respectively.
\end{definition}

\begin{proposition}
  \label{pro:quasi_isometry}
  If \(X\) and~\(Y\) are quasi\brd{}isometric, then the chain complexes
  \(\T\tilde{C}_\bullet(X)\) and \(\T\tilde{C}_\bullet(Y)\) are bornologically
  homotopy equivalent.  (This means that the chain maps and homotopies
  involved are bounded.)
\end{proposition}

\begin{proof}
  The assertion follows from the following two functoriality assertions.  Let
  \(X\) and~\(Y\) be two proper discrete metric spaces and let \(f\colon X\to
  Y\) be a quasi\brd{}Lipschitz map.  Then the induced chain map \(f_*\colon
  \tilde{C}_\bullet(X)\to\tilde{C}_\bullet(Y)\) extends to a bounded linear
  map \(f_*\colon \T\tilde{C}_\bullet(X)\to \T\tilde{C}_\bullet(Y)\).  Let
  \(f,f'\colon X\to Y\) be two quasi\brd{}Lipschitz maps that are close in the
  sense that \(d(f(x),f'(x))\) remains bounded.  Then \(H(f,f')\) extends to a
  bounded operator \(H(f,f')\colon \T\tilde{C}_\bullet(X)\to
  \T\tilde{C}_\bullet(Y)\), which still
  satisfies~\eqref{eq:Hff_chain_homotopy}.  That is, the maps \(f_*,
  f'_*\colon \T\tilde{C}_\bullet(X)\to\T\tilde{C}_\bullet(Y)\) are
  bornologically homotopy equivalent.
\end{proof}

A constant map \(X\to\pt\to X\), which is of course Lipschitz, induces the
zero map on \(\T\tilde{C}_\bullet(X)\).  However, it is not close to the
identity map unless~\(X\) is finite.  Hence \(\T\tilde{C}_\bullet(X)\) has a
chance to be non-trivial.

\section{Combings and contractibility of tempered free chain complexes}
\label{sec:combings_contract}

In this section we find a sufficient condition for \(\T\tilde{C}_\bullet(X)\)
to be contractible, based on the notion of a combing of a metric space.  This
notion is due to William P.\ Thurston and grew out of the theory of automatic
groups developed in~\cite{Epstein:Automatic}.  More recently, geometric group
theorists have become interested in more general classes of combings than
those allowed in~\cite{Epstein:Automatic} in order to overcome the following
problem: the only nilpotent groups that are combable in the sense
of~\cite{Epstein:Automatic} are the virtually Abelian ones.  It is possible to
find combings with somewhat weaker properties on certain nilpotent groups (see
\cites{Bridson:Combings_semidirect, Gilman-Holt-Rees}).  However, these
combings are asynchronous, while the argument below needs synchronous
combings.

\begin{definition}
  \label{def:combing}
  Let~\(X\) be a discrete proper metric space, choose \(\pt\in X\).  A
  \emph{combing} on~\(X\) is a sequence of maps \(f_n\colon X\to X\),
  \(n\in\N\), such that
  \begin{itemize}
  \item for all \(x\in X\), we have \(f_0(x)=\pt\) and there is \(n\in\N\)
    such that \(f_N(x)=x\) for all \(N\ge n\);
    
  \item the maps \((f_n)_{n\in\N}\) are \emph{uniformly quasi\brd{}Lipschitz},
    that is, there exists \(C>0\) with \(d(f_n(x),f_n(y))\le C\cdot
    (d(x,y)+1)\) for all \(x,y\in X\), \(n\in\N\);
    
  \item the pairs of maps \((f_n,f_{n+1})_{n\in\N}\) are \emph{uniformly
      close} in the sense that there is \(S\in\N\) such that
    \(d\bigl(f_n(x),f_{n+1}(x)\bigr)\le S\) for all \(x\in X\), \(n\in\N\).
  \end{itemize}
  
  Let \(J(x)\) be the number of \(n\in\N\) with \(f_n(x)\neq f_{n+1}(x)\).
  Since the sequence \(\bigl(f_n(x)\bigr)\) is eventually constant,
  \(J(x)<\infty\) for all \(x\in X\).  A combing has \emph{polynomial growth
    of order \(m\in\N\)} if \(J(x) \le C(\ell(x)+1)^m\) for some \(C\in\R_+\).
  \emph{Subexponential growth} and \emph{exponential growth} are defined
  similarly.
\end{definition}

The function~\(J\) is related to the length of a combing in
\cites{Gilman-Holt-Rees}.  If~\(G\) is a finitely generated discrete group,
then a combing in the above sense exists on~\(G\) if and only if~\(G\) has a
synchronous combing in the notation of \cite{Gilman-Holt-Rees}.  The
additional conditions that are imposed on combings in~\cite{Epstein:Automatic}
imply linear growth.  That is, all assertions of
Theorem~\ref{the:combable_contract} below apply to groups that are combable in
the sense of~\cite{Epstein:Automatic}.  This class of groups contains all
automatic groups and therefore all hyperbolic groups.  Thus it contains
finitely generated free non-Abelian groups and free groups.  The only
nilpotent groups that are combable in the sense of~\cite{Epstein:Automatic}
are the Abelian ones.

One checks easily that the existence of a combing on~\(X\) (with specified
growth) is independent of the choice of~\(\pt\) and a quasi\brd{}isometry
invariant.

\begin{example}
  \label{exa:free_group_combing}
  Finitely generated free groups have a combing of linear growth.
  Let~\(\F_r\) be the free group on~\(r\) generators, which we denote by
  \(s(1),\dotsc,s(r)\).  Also let \(s(r+j)=s(j)^{-1}\) for \(j=1,\dotsc,r\).
  Let~\(\ell\) be the word length for this set of generators, and let
  \(d(g,h)=\ell(g^{-1}h)\) be the associated left invariant distance.  We
  define the sequence of maps \((f_n)\) as follows.  Write \(g\in\F_r\) as a
  reduced word \(g=s(i_1)\dotsm s(i_\ell)\) in the generators, that is,
  \(\ell=\ell(g)\).  Let \(f_n(g)=s(i_1)\dotsm s(i_n)\) for \(0\le n<\ell\)
  and \(f_n(g)=g\) for \(n\ge\ell\).  Thus \(J(g)=\ell(g)\) for this sequence
  of maps.  One checks easily that \(d(f_n(g),f_{n+1}(g))\le 1\) for all
  \(g\in \F_r\), \(n\in\N\), and \(d(f_n(g),f_n(h))\le d(g,h)\) for all
  \(g,h\in\F_r\), \(n\in\N\).  Hence we have a combing of linear growth.
\end{example}

\begin{theorem}
  \label{the:combable_contract}
  Let~\(X\) be a discrete proper metric space.  If~\(X\) has a combing of
  polynomial growth of order \(m\in\N\), then the embeddings
  \(\Sch^{k+m}\tilde{C}_\bullet(X)\to\Sch^k\tilde{C}_\bullet(X)\) are
  homotopic to zero for all \(k\in\R_+\).
  
  If~\(X\) has a combing of polynomial growth, then
  \(\Sch^\infty\tilde{C}_\bullet(X)\) is contractible.
  
  If~\(X\) has a combing of subexponential growth, then
  \(\Sch^\omega\tilde{C}_\bullet(X)\) is contractible.
  
  If~\(X\) has a combing of exponential growth, then
  \(\Hol\tilde{C}_\bullet(X)\) is contractible.
\end{theorem}

\begin{proof}
  We only prove the first assertion; the others are proven similarly.  The
  definition of \(H(f,f')\) shows that \(H(f,f')(x_0,\dotsc,x_n)=0\) unless
  \(f(x_j)\neq f'(x_j)\) for some \(j\in\{0,\dotsc,n\}\).  Therefore, for a
  basis vector \((x_0,\dotsc,x_n)\) there are only \(J(x_0)+\dotsb+J(x_n)\)
  non-zero terms in the sum
  \[
  H(x_0,\dotsc,x_n) \defeq \sum_{i=0}^\infty H(f_i,f_{i+1})(x_0,\dotsc,x_n).
  \]
  The linear map \(H\colon \tilde{C}_\bullet(X)\to \tilde{C}_\bullet(X)\) of
  degree~\(1\) so defined satisfies
  \[
  \delta H+ H \delta
  = \sum_{i=0}^\infty (f_{i+1})_* - (f_i)_*
  = \lim_{i\to\infty} (f_i)_* - (f_0)_*
  = \ID
  \]
  by~\eqref{eq:Hff_chain_homotopy}.
  That is, it is a contracting homotopy for \(\tilde{C}_\bullet(X)\).
  If~\(H\) extends to a bounded linear map \(H\colon
  \Sch^{k+m}\tilde{C}_\bullet(X)\to \Sch^k\tilde{C}_\bullet(X)\), then it is
  the desired chain homotopy between the zero map and the embedding
  \(\Sch^{k+m}\tilde{C}_\bullet(X)\to \Sch^k\tilde{C}_\bullet(X)\).
  
  Since the maps \(f_i, f_{i+1}\) are uniformly quasi\brd{}Lipschitz and
  close, the control condition on elements of \(\Sch^k\tilde{C}(X)\) is
  uniformly preserved by all operators \(H(f_i,f_{i+1})\) and hence by~\(H\).
  It therefore remains to consider the growth condition.  For any \(R>0\), the
  unit ball of \(\Sch^{k+m} C_n(X)_R\) is the complete absolutely convex hull
  of the rescaled basis vectors \((\ell(x_0)+1)^{-k-m}\cdot (x_0,\dots,x_n)\)
  with \(x_0,\dots,x_n\in X\) satisfying \(d(x_i,x_j)\le R\).  Hence~\(H\) is
  bounded on \(\Sch^{k+m} C_n(X)_R\) if (and only if) the images of these
  rescaled basis vectors form a uniformly bounded subset of \(\Sch^k
  C_{n+1}(X)\).
  
  There is \(C>0\) with \(\ell(f_i(x))+1\le C(\ell(x)+1)\) for all \(i\in\N\),
  \(x\in X\).  This follows from the uniform quasi\brd{}Lipschitz condition
  and the finiteness of the set of \(f_i(*)\), \(i\in\N\).  Hence the set of
  \[
  (-1)^j (\ell(x_0)+1)^{-k}\cdot
  (f_i(x_0),\dotsc,f_i(x_j),f_{i+1}(x_j),\dotsc,f_{i+1}(x_n))
  \]
  for \(x_0,\dotsc,x_n\in X\) with \(d(x_i,x_j)\le R\) is uniformly bounded in
  \(\Sch^k C_n(X)\).  The image of \((\ell(x_0)+1)^{-k}\cdot
  H(x_0,\dots,x_n)\) is a sum of \(J(x_0)+\dotsb+J(x_n)\) of these basis
  vectors.  This number is \(O(\ell(x_0)+1)^m\) by hypothesis.  Hence the set
  of \((\ell(x_0)+1)^{-k-m}\cdot H(x_0,\dots,x_n)\) for \(x_0,\dotsc,x_n\in
  X\) with \(d(x_i,x_j)\le R\) is uniformly bounded in \(\Sch^k C_{n+1}(X)\)
  as desired.
\end{proof}

\section{Tempered group cohomology and projective resolutions}
\label{sec:tempered_cohomology}

Let~\(G\) be a finitely generated discrete group.  Recall the following
standard definition of \(H^n(G)\).  Let \(\Coh^n(G)\) be the space of all
functions \(G^n\to\R\), with the convention \(\Coh^0(G)=\R\); define
\(\delta\colon \Coh^n(G)\to\Coh^{n+1}(G)\) by
\begin{multline*}
  \delta\phi(g_1,\dotsc,g_{n+1}) \defeq \phi(g_2,\dotsc,g_n)
  \\ + \sum_{j=1}^{n} (-1)^j \phi(g_1,\dotsc,g_j\cdot g_{j+1},\dotsc,g_{n+1})
  + (-1)^{n+1} \phi(g_1,\dotsc,g_n).
\end{multline*}
Then \(\Coh^\bullet(G)\defeq (\Coh^n(G),\delta)_{n\ge0}\) is a chain complex
with \(H^n(\Coh^\bullet(G))\cong H^n(G)\).

We are interested in subcomplexes of \(\Coh^\bullet(G)\) defined by growth
conditions.  Let \(\Coh^n(\Sch^k G) \subseteq \Coh^n(G)\) be the subspace of
cochains \(\phi\colon G^n\to\R\) that satisfy
\[
\abs{\phi(g_1,\dotsc,g_n)}\le C (\ell(g_1)+1)^k \dotsm (\ell(g_n)+1)^k
\]
for some \(C>0\) and for \(k\in\R_+\).  We let \(\Coh^n(\Sch^\infty G) =
\bigcup_{k\in\N} \Coh^n(\Sch^k G)\).  We let \(\Coh^n(\Sch^\omega G)\) be the
subspace of \(\phi\colon G^n\to\R\) that satisfy
\begin{equation}
  \label{eq:exponential_growth}
  \abs{\phi(g_1,\dotsc,g_n)}\le C \alpha^{\ell(g_1)} \dotsm \alpha^{\ell(g_n)}
  = C \alpha^{\ell(g_1)+\dotsb+\ell(g_n)}
\end{equation}
for all \(\alpha>1\).  We define \(\Coh^n(\Hol G)\) to consist of those
\(\phi\colon G^n\to\R\) that satisfy~\eqref{eq:exponential_growth} for some
\(\alpha>1\).  One checks immediately that the boundary map~\(\delta\)
preserves these growth conditions, so that we get chain subcomplexes of
\(\Coh^\bullet(G)\).  We let \(H^*(\Sch^k G)\) for
\(k\in\R_+\cup\{\infty,\omega\}\), \(*\in\N\), and \(H^*(\Hol G)\) for
\(*\in\N\) be the resulting cohomology groups.  The cohomology \(H^*(\Sch^0
G)=H^*(\ell_1 G)\) is the \emph{bounded cohomology} of~\(G\).  We call
\(H^*(\Sch^\infty G)\) \emph{polynomial growth cohomology}, \(H^*(\Sch^\omega
G)\) \emph{subexponential growth cohomology}, and \(H^*(\Hol G)\)
\emph{exponential growth cohomology}.  There are natural maps \(H^*(\Sch^k G)
\to H^*(\Sch^l G) \to H^*(\Hol G)\to H^*(G)\) for
\(k,l\in\R_+\cup\{\infty,\omega\}\) with \(l\ge k\).  We want to find
sufficient conditions for some of them to be isomorphisms.

First we rewrite the chain complexes above in terms of bar resolutions.
Let~\(A\) be a bornological unital algebra and let~\(M\) be a bornological
left module over~\(A\).  The \emph{bar resolution} \(\Barr_\bullet(A,M)\)
of~\(M\) is defined as follows.  Let \(\Barr_n(A,M)\defeq A^{\hot n+1} \hot
M\) for \(n\ge0\) and view this as a left \(A\)\nbd{}module by \(a\cdot
(x_0\otimes \dotsb \otimes x_{n+1}) \defeq (a\cdot x_0)\otimes x_1\otimes
\dotsb\otimes x_{n+1}\).  Here~\(\hot\) denotes the projective complete
bornological tensor product (see~\cite{Hogbe-Nlend:Completions}).  We also
define \(\Barr_{-1}(A,M)=M\).  The boundary map \(b'\colon \Barr_n(A,M) \to
\Barr_{n-1}(A,M)\) is defined by
\[
b'(x_0\otimes\dotsb\otimes x_{n+1}) \defeq
\sum_{j=0}^n (-1)^j x_0\otimes\dotsb\otimes x_j\cdot x_{j+1}
\otimes\dotsb\otimes x_{n+1}
\]
for all \(n\ge0\), \(x_0,\dots,x_n\in A\), \(x_{n+1}\in M\).  Then
\((b')^2=0\) and~\(b'\) is a left module homomorphism, so that we have a chain
complex of bornological left \(A\)\nbd{}modules.  The map
\(s(x_0\otimes\dotsb\otimes x_n)\defeq 1\otimes x_0\otimes\dotsb\otimes x_n\)
is a bounded contracting homotopy, that is, \(s b'+b's = \ID\).  Moreover,
\(\Barr_n(A,M)\) is a free module for all \(n\ge0\).  Thus the bar resolution
is a free bornological left \(A\)\nbd{}module resolution of~\(M\).  Let~\(N\)
be another bornological left \(A\)\nbd{}module and let
\(\Hom_A(\Barr_\bullet(A,M),N)\) be the chain complex of bounded left module
homomorphisms \(f\colon \Barr_n(A,M)\to N\) for \(n\in\N\).  We let
\(\Ext^n_A(M,N)\) be the \(n\)th cohomology of
\(\Hom_A(\Barr_\bullet(A,M),N)\).

Since
\[
\Hom_A (A\hot X, N)\cong \Hom(X,N)
\]
for any bornological vector space~\(X\) and any bornological left
\(A\)\nbd{}module~\(N\), free modules are projective for extensions with a
bounded linear section.  The same proof as in Abelian categories shows that
any two free resolutions of the same module are homotopy equivalent via
bounded \(A\)\nbd{}linear chain maps and homotopies.  Therefore, if
\(P_\bullet\to M\) is any free resolution of~\(M\) (this means that we have a
bounded contracting homotopy), then the chain complex \(\Hom_A(P_\bullet,N)\)
is naturally homotopy equivalent to \(\Hom_A(\Barr_\bullet(A,M),N)\), so that
\(\Ext^n_A(M,N)\) is naturally isomorphic to the cohomology of the chain
complex \(\Hom_A(P_\bullet,N)\).

\begin{lemma}
  \label{lem:tempered_group_cohomology_Ext}
  Equip~\(\R\) with the trivial representation of~\(G\).  Then
  \[
  H^n(G) \cong \Ext^n_{\R[G]}(\R,\R),
  \qquad
  H^n(\T G) \cong \Ext^n_{\T(G)}(\R,\R)
  \]
  for \(\T(G)=\Hol(G)\) or \(\T(G)=\Sch^k(G)\),
  \(k\in\R_+\cup\{\infty,\omega\}\).
\end{lemma}

\begin{proof}
  We write \(H^n(\R[G]) = H^n(G)\) and allow \(\T(G)=\R[G]\) among our list of
  convolution algebras.  We claim that the chain complex
  \(\Hom_{\T(G)}(\Barr_\bullet(\T(G),\R),\R)\) is canonically isomorphic to
  the chain complexes that we have used above to define \(H^n(\T G)\).  This
  establishes the lemma.
  
  Consider \(\R[G]\) first.  Then \(\Barr_n(\R[G],\R)=\R[G]^{\hot n+1}\).
  This is isomorphic to \(\R[G^{n+1}]\) with the fine bornology.  A
  \(G\)\nbd{}equivariant linear map to~\(\R\) is the same as a
  \(G\)\nbd{}invariant linear functional.  Such linear functionals on
  \(\R[G^{n+1}]\) are of the form
  \begin{equation}
    \label{eq:function_to_functional}
    \psi \mapsto
    \sum_{g_0,\dotsc,g_n\in G} \psi(g_0,\dotsc,g_n) \phi(g_1,\dotsc,g_n)
  \end{equation}
  for a uniquely determined function \(\phi\colon G^n\to\R\).  Thus
  \(\Hom_{\R[G]}(\Barr_n(\R[G],\R),\R)\cong \Coh^n(G)\).  One checks
  immediately that these isomorphisms form a chain map.
  
  For another choice of \(\T(G)\), the complex \(\Barr_\bullet(\T(G),\R)\)
  contains \(\Barr_\bullet(\R[G],\R)\) as a dense subcomplex.  Hence we may
  view \(\Hom_{\T(G)}(\Barr_\bullet(\T(G),\R),\R)\) as a subcomplex of
  \(\Hom_{\R[G]}(\Barr_\bullet(\R[G],\R),\R) \cong \Coh^\bullet(G)\).  We
  claim that it agrees with the subcomplex that defines \(H^n(\T G)\).  The
  space \(\Hom_{\T(G)}(\Barr_n(\T(G),\R),\R)\) is the space of all
  \(\phi\colon G^n\to\R\) for which~\eqref{eq:function_to_functional} defines
  a bounded linear functional on \(\T(G)^{\hot n}\).  Equivalently, the
  \(n\)\nbd{}linear map
  \[
  (\psi_1,\dots,\psi_n)\mapsto \sum_{g_1,\dots,g_n\in G} \psi_1(g_1)\dotsm
  \psi_n(g_n) \phi(g_1,\dotsc,g_n)
  \]
  is bounded on \(\T(G)^n=\T(G)\times \dotsb\times \T(G)\).  We have identified the dual
  space of \(\T(G)\) in Section~\ref{sec:function_spaces}.  The same
  discussion applies to multilinear maps and yields the assertions of the
  lemma.
\end{proof}

\begin{theorem}
  \label{the:tempered_cohomology}
  Let~\(G\) be a finitely generated discrete group.  If
  \(\T\tilde{C}_\bullet(G)\) has a bounded contracting homotopy, then the
  canonical map \(H^n(\T G)\to H^n(G)\) is an isomorphism.
\end{theorem}

\begin{proof}
  Let~\(G\) act diagonally on \(\T C_n(G)\) by \(g\cdot (x_0,\dotsc,x_n)\defeq
  (gx_0,\dotsc,gx_n)\).  Let \(C_n'(G)\subseteq C_n(G)\) be the subspace
  spanned by \((1,x_1,\dotsc,x_n)\).  Let \(C_n'(G)_R\) for \({R>0}\) be the
  finite-dimensional subspace spanned by \((1,x_1,\dotsc,x_n)\) with
  \(d(x_i,x_j)\le R\) and \(d(1,x_i)\le R\) for all \(i,j\in\{1,\dots,n\}\).
  There is a bornological isomorphism
  \[
  \T(G) \hot C_n'(G)_R \to \T C_n(G)_R,
  \qquad
  g\otimes (1,x_1,\dotsc,x_n) \mapsto (g,gx_1,\dotsc,gx_n)
  \]
  because of the left invariance of the metric on~\(G\).  This implies
  \[
  \T C_n(G)
  \cong \varinjlim \T(G) \hot C_n'(G)_R
  \cong \T(G) \hot C_n'(G)
  \]
  because~\(\hot\) commutes with direct limits.  If we give \(\T(G)\hot
  C_n'(G)\) the standard \(\T(G)\)\brd{}module structure via \(\phi\cdot
  (\phi'\otimes \psi) \defeq (\phi*\phi')\otimes\psi\), then the above
  isomorphism is a module homomorphism.  Hence \(\T C_n(G)\) is a free
  \(\T(G)\)\brd{}module.  It is easy to see that the boundary maps of \(\T
  C_\bullet(G)\) and the augmentation map \(\T C_0(G)\to \R\) are
  \(G\)\nbd{}equivariant and hence \(\T(G)\)\brd{}module homomorphisms.
  
  Since \(\tilde{C}_\bullet(G)\) is contractible, \(C_\bullet(G)\to\R\) is a
  free \(\R[G]\)\brd{}module resolution of~\(\R\).  (It is isomorphic to the
  \emph{reduced} bar resolution.)  Since all free resolutions are equivalent,
  Lemma~\ref{lem:tempered_group_cohomology_Ext} yields \(H^n(G) \cong
  H^n(\Hom_{\R[G]}(C_\bullet(G),\R))\).  The canonical isomorphisms
  \[
  \Hom_{\R[G]}(\R[G]\hot C_n'(G),\R) \cong
  \Hom(C_n'(G),\R) \cong
  \Hom_{\T(G)}(\T(G)\hot C_n'(G),\R)
  \]
  are compatible with the boundary maps on \(C_\bullet(G)\) and \(\T
  C_\bullet(G)\), that is, they form an isomorphism of chain complexes
  \begin{equation}
    \label{eq:compare_chain_complexes}
    \Hom_{\R[G]}(C_\bullet(G),\R) \cong \Hom_{\T(G)}(\T C_\bullet(G),\R).
  \end{equation}
  If \(\T\tilde{C}_\bullet(G)\) has a bounded contracting homotopy, then \(\T
  C_\bullet(G)\) augmented by \(\alpha\colon \T C_0\to\R\) is a free
  \(\T(G)\)\brd{}module resolution of~\(\R\).  By the uniqueness of free
  resolutions, this implies \(H^n(\T G) \cong H^n(\Hom_{\T(G)}(\T
  C_\bullet(G),\R))\).  The assertion of the theorem now follows
  from~\eqref{eq:compare_chain_complexes}.
\end{proof}

\begin{corollary}
  \label{cor:tempered_cohohomology_combable}
  If~\(G\) has a combing of polynomial growth, then \(H^n(\Sch^\infty(G))
  \cong H^n(G)\).  If~\(G\) has a combing of subexponential growth, then
  \(H^n(\Sch^\omega(G)) \cong H^n(G)\).  If~\(G\) has a combing of exponential
  growth, then \(H^n(\Hol(G)) \cong H^n(G)\).
\end{corollary}

\begin{proof}
  Combine Theorems \ref{the:combable_contract}
  and~\ref{the:tempered_cohomology}.
\end{proof}

The following theorem provides an obstruction for the contractibility of
\(\T\tilde{C}_\bullet(G)\).  Actually, the proof shows that it even yields an
obstruction to \(H^n(\T G)\to H^n(G)\) being a bornological isomorphism.

\begin{theorem}
  \label{the:finite_cohomology}
  Let~\(G\) be a finitely generated discrete group and let \(\T(G)\) be one of
  \(\Sch^\infty(G)\), \(\Sch^\omega(G)\), or \(\Hol(G)\).  Suppose that
  \(\T\tilde{C}_\bullet(G)\) is contractible.  Then \(H^n(G)\) is
  finite-dimensional for all \(n\in\N\).
\end{theorem}

\begin{proof}
  By Theorem~\ref{the:tempered_cohomology}, \(H^n(G)\) and \(H^n(\T G)\) are
  isomorphic.  Even more, they are computed by chain homotopy equivalent
  bornological chain complexes.  The cohomology of a bornological chain
  complex carries a canonical bornology, which may fail to be separated.  The
  homotopy equivalence of the underlying bornological chain complexes implies
  that \(H^n(G)\) and \(H^n(\T G)\) are isomorphic as bornological vector
  spaces.  Replacing \(\T(G)\) by \(\Hol(G)\), we no longer necessarily have
  an isomorphism, but we obtain that \(H^n(G)\) is a direct summand in
  \(H^n(\Hol G)\) because \(\R[G]\subseteq \Hol(G)\subseteq \T(G)\).
  
  The group cohomology \(H^n(G)\) is the dual space of the group homology
  of~\(G\).  The latter is isomorphic to \(\R[B]\) for some countable
  set~\(B\).  The dual space of \(\R[B]\) is \(\prod_{x\in B} \R\), equipped
  with the product bornology.  We get a bornological isomorphism \(H^n(G)\cong
  \prod_{x\in B} \R\).  Thus \(H^n(G)\) is a Fréchet-Schwartz space.  In
  contrast, the complex that computes \(H^n(\Hol G)\) is a complex of Silva
  spaces.  Hence any separated subspace of \(H^n(\Hol G)\) is a Silva space.
  In particular, \(H^n(G)\) is a Fréchet-Schwartz and a Silva space at the
  same time.  We claim that such a space is necessarily finite-dimensional.
  
  Let~\(V\) be a Fréchet-Schwartz space whose von Neumann bornology is a Silva
  bornology.  Thus \(V=\varinjlim V_n\) for an inductive system
  \((V_n)_{n\in\N}\) of Banach subspaces with injective maps \(V_n\to
  V_{n+1}\).  Let~\(S_n\) be the unit ball of~\(V_n\).  Since~\(V\) is a
  Fréchet space, there exists a sequence of positive scalars \((\epsilon_n)\)
  such that \(\bigcup \epsilon_n S_n\) is still von Neumann bounded.  This is
  the bornological definition of metrisability, see~\cite{Meyer:Born_Top}.
  However, any bounded subset of \(\varinjlim V_n\) is absorbed by~\(S_N\) for
  some \({N\in\N}\).  Hence~\(S_N\) absorbs~\(S_n\) for all \(n\in\N\).
  Thus~\(V\) is isomorphic to the Banach space~\(V_N\).  Since~\(V\) is
  Fréchet-Schwartz, any bounded subset in~\(V\) is relatively compact.
  However, the only Banach spaces with relatively compact unit ball are the
  finite-dimensional ones.  Hence~\(V\) is finite-dimensional as asserted.
\end{proof}

It is known that \(BG\) has \emph{finite type} if~\(G\) is combable.  That is,
\(BG\) has the homotopy type of a simplicial complex with finite skeleta.  Of
course, this implies that \(H^n(G)\cong H^n(BG)\) is finite-dimensional.
Other popular conditions of geometric group theory like amenability or
exactness do not guarantee \(BG\) to have finite-dimensional cohomology spaces
and hence cannot possibly suffice to prove contractibility of
\(\T\tilde{C}_\bullet(G)\).  Only conditions that guarantee some finiteness of
\(BG\) have a chance to imply contractibility of \(\T\tilde{C}_\bullet(G)\).

\section{Using other free resolutions and the example of free groups}
\label{sec:other_resolution}

Let~\(G\) be a finitely generated discrete group and let \(\T(G)\) be one of
the algebras \(\Hol(G)\) or \(\Sch^k(G)\) for
\(k\in\R_+\cup\{\infty,\omega\}\).  We want to find simpler chain complexes
that are homotopy equivalent to \(\T C_\bullet(G)\).

Recall that \(\Mod(A)\) denotes the category of bornological left
\(A\)\nbd{}modules.  The embedding \(i\colon \R[G]\to\T(G)\) induces a
``forgetful'' functor
\[
i^*\colon \Mod(\T(G)) \to \Mod(\R[G]).
\]
We claim that~\(i^*\) has a left adjoint functor \(i_!\colon \Mod(\R[G]) \to
\Mod(\T(G))\).  That is, the functors \(i^*\) and~\(i_!\) are related by
natural isomorphisms
\begin{equation}
  \label{eq:induced_adjoint}
  \Hom_{\T(G)}(i_!(M),N) \cong \Hom_{\R[G]}(M,i^*(N))
\end{equation}
if \(M\) and~\(N\) are bornological modules over \(\R[G]\) and \(\T(G)\),
respectively.  We construct~\(i_!\) as follows (see
also~\cite{Meyer:Embed_derived}).  Let~\(M\) be a left bornological
\(\R[G]\)\brd{}module.  We let \(\T(G)\hot_G M\) be the quotient of
\(\T(G)\hot M\) by the closed linear subspace generated by
\(\phi\ast\delta_g\otimes m - \phi\otimes gm\) for \(\phi\in\T(G)\), \(g\in
G\), \(m\in M\).  Since we take closures, this quotient is again a complete
convex bornological vector space.  We turn \(\T(G)\hot_G M\) into a left
bornological \(\T(G)\)\brd{}module by \(\phi_1\cdot (\phi_2\otimes m)\defeq
\phi_1\ast\phi_2\otimes m\) for \(\phi_1,\phi_2\in\T(G)\), \(m\in M\).  We
claim that \(i_!(M)=\T(G)\hot_G M\).

To prove this, we verify~\eqref{eq:induced_adjoint}.  Let~\(N\) be a
bornological \(\T(G)\)\brd{}module.  By the universal property of~\(\hot\),
bounded \(\T(G)\)\brd{}module homomorphisms \(\T(G)\hot_G M\to N\) correspond
to bilinear maps \(f\colon \T(G)\times M\to N\) that satisfy
\[
f(\phi_1\ast\phi_2, m) = \phi_1\cdot f(\phi_2, m),
\qquad
f(\phi\ast\delta_g, m) = f(\phi, g\cdot m)
\]
for all \(\phi,\phi_1,\phi_2\in\T(G)\), \(g\in G\), \(m\in M\).  The first
relation implies that \(f(\phi, m) = \phi\cdot f(\delta_1, m)\), the second
relation implies that the bounded linear map \(M\to N\), \(m\mapsto
f(\delta_1, m)\), is \(G\)\nbd{}equivariant.  Conversely, given any bounded,
\(G\)\nbd{}equivariant linear map \(f'\colon M\to N\), the formula \(f(\phi,
m)\defeq \phi\cdot f'(m)\) defines a bounded \(\T(G)\)\brd{}module
homomorphism \(\T(G)\hot_G M\to N\).  Thus we have constructed the natural
isomorphism~\eqref{eq:induced_adjoint}.  We usually write \(\T(G)\hot_G M\)
and not \(i_!(M)\) because we want to keep \(\T(G)\) in our notation.

\begin{proposition}
  \label{pro:chains_resolutions}
  There is a canonical isomorphism of bornological chain complexes
  \[
  \T(G) \hot_G C_\bullet(G) \cong \T C_\bullet(G).
  \]
  If \(P_\bullet\to\R\) is any free \(\R[G]\)\brd{}module resolution of~\(G\),
  then \(\T(G) \hot_G P_\bullet\) is bornologically homotopy equivalent to
  \(\T C_\bullet(G)\) (the chain maps and chain homotopies involved are
  bounded linear maps).
\end{proposition}

\begin{proof}
  During the proof of Theorem~\ref{the:tempered_cohomology}, we have shown
  \[
  C_n(G) \cong \R[G]\hot C_n'(G),
  \qquad
  \T C_n(G) \cong \T(G)\hot C_n'(G),
  \]
  for certain fine bornological vector spaces \(C_n'(G)\); these are module
  isomorphisms with respect to the free module structures on the right hand
  sides.  It follows easily that \(\T C_n(G) \cong \T(G)\hot_G C_n(G)\).
  These isomorphisms form an isomorphism of chain complexes \(\T C_\bullet(G)
  \cong \T(G)\hot_G C_\bullet(G)\).  This implies the second assertion because
  of the uniqueness up to homotopy equivalence of free resolutions.
\end{proof}

\begin{example}
  \label{exa:free_group_resolution}
  Let \(r\in\N_{\ge1}\) and let~\(\F_r\) be the non-Abelian free group
  on~\(r\) generators, denoted \(s(1),\dotsc,s(r)\).  Let \(\T(\F_r)\) be one
  of the convolution algebras considered above.
  Example~\ref{exa:free_group_combing} provides a combing of linear growth
  on~\(\F_r\).  Theorem~\ref{the:combable_contract} yields that
  \(\T\tilde{C}_\bullet(\F_r)\) is contractible if \(\T(\F_r)\) is
  \(\Sch^\infty(\F_r)\), \(\Sch^\omega(\F_r)\), or \(\Hol(\F_r)\).  We are
  going to show that this fails for \(\Sch^k(\F_r)\) with \(k\in\R_+\) and,
  especially, for \(\ell_1(\F_r)\).
  
  The Cayley graph of~\(\F_r\) with respect to the generators above is a tree.
  The simplicial chain complex associated to this tree is a free
  \(\R[\F_r]\)\brd{}module resolution of the trivial representation.
  Explicitly, this yields the free resolution
  \begin{equation}
    \label{eq:free_group_resolution}
    \dotsb \longrightarrow 0 \longrightarrow \R[\F_r]^r
    \overset{b_1}\longrightarrow \R[\F_r]
    \overset{b_0}\longrightarrow  \R \longrightarrow 0,
  \end{equation}
  where
  \[
  b_0(\phi) = \sum_{g\in\F_r} \phi(g),
  \qquad
  b_1(\phi_1,\dotsc,\phi_r)(g)
  = \sum_{j=1}^r \phi_j(g) - \phi_j\bigl(gs(j)\bigr).
  \]
  It is clear that the maps \(b_0\) and~\(b_1\) are module homomorphisms and
  that \(b_0b_1=0\).  To check by hand that~\eqref{eq:free_group_resolution}
  is a resolution, we specify a linear map \(\sigma\colon \R[\F_r]\to
  \R[\F_r]^r\) that satisfies \(\sigma b_1=\ID\) and \(b_1\sigma(g) = g-1\)
  for all \(g\in \F_r\).  Actually, such a map is necessarily unique.  We can
  obtain \(\sigma(g)\) for \(g\in\F_r\) by induction on the length of~\(g\).
  If \(g=1\), then \(\sigma(g)=0\).  Otherwise, write \(g=g_2s(j)^\epsilon\)
  with \(j\in\{1,\dotsc,r\}\), \(\epsilon\in\{\pm 1\}\), and
  \(\ell(g_2)<\ell(g)\).  Then we let \(\sigma(g)-\sigma(g_2)\) be equal
  to~\(g\) in the \(j\)th summand if \(\epsilon=1\), and~\(-g_2\) in the
  \(j\)th summand if \(\epsilon=-1\).  This determines \(\sigma(g)\) for all
  \(g\in\F_r\).  The relation \(\sigma b_1=\ID\) is trivial, and
  \(b_1\sigma(g)=g-1\) follows by induction on \(\ell(g)\).
  
  Proposition~\ref{pro:chains_resolutions} yields that \(\T C_\bullet(\F_r)\)
  is chain homotopy equivalent to
  \begin{equation}
    \label{eq:free_group_chains}
    \dotsb \longrightarrow 0 \longrightarrow \T(\F_r)^r
    \overset{b_1}\longrightarrow \T(\F_r).
  \end{equation}
  A chain complex equivalent to \(\T\tilde{C}_\bullet(\F_r)\) is obtained by
  replacing \(\T(\F_r)\) by the kernel of the augmentation map \(b_0\colon
  \T(\F_r)\to\R\).  As a result, \(\T\tilde{C}_\bullet(\F_r)\) is split exact
  except possibly in dimensions \(0\) and~\(1\).  It is split exact in all
  dimensions if and only if~\(b_1\) is a bornological isomorphism from
  \(\T(\F_r)^r\) onto \(\ker b_0 \subseteq \T(\F_r)\).
  
  The construction above shows that \(\sigma(g)\) is a sum of exactly
  \(\ell(g)\) different standard basis vectors of \(\R[\F_r]^r\).
  Hence~\(\sigma\) is a bounded linear map for \(\Hol(\F_r)\),
  \(\Sch^\infty(\F_r)\), and \(\Sch^\omega(\F_r)\), and maps
  \(\Sch^{k+1}(\F_r)\to\Sch^k(\F_r)^r\) for \(k\in\R_+\).  These assertions
  are already contained in Theorem~\ref{the:combable_contract}.  The
  map~\(\sigma\) is, in fact, closely related to the combing for~\(\F_r\)
  constructed in Example~\ref{exa:free_group_combing}.
  
  Let \(\T(\F_r) = \Sch^k(\F_r)\) for \(k\in\R_+\).  Since~\(\sigma\) is a
  bounded linear map from \(\Sch^k(\F_r)\) to \(\Sch^{k-1}(\F_r)^r\), the
  map~\(b_1\) is injective on \(\Sch^k(\F_r)^r\).  Therefore, the chain
  complex \(\T C_\bullet(\F_r)\) is split exact in degree~\(1\) in the sense
  that the boundary map \(\Sch^k C_2(\F_r)\to \ker \delta\subseteq \Sch^k
  C_1(\F_r)\) is split surjective.  Since~\eqref{eq:free_group_resolution} is
  exact, \(b_1(\Sch^k(\F_r)^r)\) contains \(b_1(\R[\F_r]^r)=\ker b_0\cap
  \R[\F_r]\), which is dense in \(\ker b_0\colon \Sch^k(\F_r)\to\R\).  We have
  \(b_1(\Sch^k(\F_r)^r)\neq \ker b_0\) already for \(r=1\), that is, for the
  group of integers; this amounts to a well-known fact about Fourier series.
  Thus \(H_0(\Sch^k\tilde{C}_\bullet (\F_r))\) is non-trivial, but its
  separated quotient (called reduced homology
  in~\cite{Mineyev:Homology_Combable}) vanishes.
\end{example}

\section{Tempered cellular chain complexes}
\label{sec:tempered_cellular_chains}

Proposition~\ref{pro:chains_resolutions} allows us to replace \(\T
C_\bullet(G)\) by simpler chain complexes for a discrete group~\(G\).  For
more general metric spaces, we do not have the category of
\(\R[G]\)\brd{}modules and hence cannot speak of projective resolutions.
Therefore, we have to work harder to find simpler models for \(\T
C_\bullet(X)\).  In this section, we relate \(\T C_\bullet(X)\) to tempered
variants of cellular chain complexes.  The motivation for this is two-fold.
These complexes sometimes arise in applications, and we may want to use
Theorem~\ref{the:combable_contract} to get information about them.
Conversely, if~\(X\) is not known to be combable, we may want to check the
contractibility of \(\T\tilde{C}_\bullet(X)\) for a smaller complex that is
homotopy equivalent to it.

As a preparation, we briefly recall the construction of cellular chain
complexes.  Let~\(X\) be a CW\brd{}complex.  We denote the set of closed
\(n\)\nbd{}cells in~\(X\) by \(X_n\) and write \(X^{(n)}\) for the \(n\)th
skeleton of~\(X\).  Thus \(X^{(0)}\subseteq X\) is a discrete subset.  A
continuous map \(f\colon X\to Y\) is called \emph{cellular} if it maps
\(X^{(n)}\) to \(Y^{(n)}\) for all \(n\in\N\).  Cellular homotopies are
defined using a canonical cell decomposition of \(X\times [0,1]\).  It is well-known
that any continuous map between CW\brd{}complexes is homotopic to a cellular
map and that homotopic cellular maps are homotopic by a cellular homotopy.

Let \(f\colon X\to Y\) be cellular.  For cells \(D\in X_n\), \(E\in Y_n\), one
defines the \emph{incidence number} \(\iota(f,D,E)\) using the action of~\(f\)
on the singular homology of the CW\brd{}pairs \((X^{(n)},X^{(n-1)})\) and
\((Y^{(n)},Y^{(n-1)})\).  Let \(\iota(\partial D,F)\) for \(D\in X_{n+1}\),
\(F\in X_n\) be the incidence number of the characteristic map
\(S^n\to\partial D\subseteq X\).  The \emph{cellular chain complex}
\((\Cell_\bullet(X),\delta)_{n\ge0}\) of~\(X\) is defined by
\(\Cell_n(X)=\R[X_n]\) and
\[
\delta\colon \Cell_{n+1}(X)\to\Cell_n(X),
\qquad
\delta(D) \defeq \sum_{F\in X_n} \iota(\partial D,F) F.
\]
A cellular map \(f\colon X\to Y\) induces a chain map
\[
f_*\colon \Cell_\bullet(X) \to \Cell_\bullet(Y),
\qquad
f_*(D) \defeq \sum_{E\in Y_n} \iota(f,D,E) E.
\]
This construction is functorial, that is, \(\ID_*=\ID\) and \((fg)_*=f_*g_*\).
A cellular homotopy \(F\colon X\times [0,1]\to Y\) induces a chain homotopy
\(H(F)\colon \Cell_\bullet(X)\to\Cell_\bullet(Y)\) between \((F_0)_*\) and
\((F_1)_*\), that is, \(\delta H(F)+H(F)\delta = (F_1)_*-(F_0)_*\).

For a cellular map \(f\colon X\to Y\) and \(n\in\N\), we let
\begin{equation}
  \label{eq:cellular_norm}
  \norm{f}_n \defeq \sup_{D\in X_n} \sum_{E\in Y_n} \abs{\iota(f,D,E)}
\end{equation}
We define similar norms for cellular homotopies and the boundary
map~\(\delta\).  We call \(f\), \(F\), or~\(\delta\) \emph{bounded} if all
these norms remain finite.  For infinite CW\brd{}complexes, this depends on
the chosen cell decomposition: if \((D_j)\) is a sequence of different
\(n\)\nbd{}cells and if we subdivide~\(D_j\) into \(j\)~cells for each~\(j\),
then the identity map becomes an unbounded cellular map.  Boundedness
of~\(\delta\) is automatic for simplicial and cubical cell complexes because
the boundaries of an \(n\)\nbd{}simplex or an \(n\)\nbd{}cube have at most
\(n+1\) or \(2^n\) nondegenerate faces, respectively.

Recall that the chain complex \(C_\bullet(X)\) for a set~\(X\) is the reduced
simplicial chain complex of a certain simplicial set \(S(X)\).  If we view the
simplicial realisation \(\abs{S(X)}\) as a CW\brd{}complex, then we get
\(C_\bullet(X) = \Cell_\bullet(\abs{S(X)})\).  We now describe some additional
structure on \(\abs{S(X)}\) that allows us to define the tempered chain
complexes \(\T C_\bullet(X)\).  First, we use the filtration \(\abs{S(X)} =
\bigcup \abs{S(X)_R}\), where \(S(X)_R\) is the simplicial subset whose
nondegenerate \(n\)\nbd{}simplices are only those \((x_0,\dots,x_n)\) with
\(n\le R\) and \(d(x_i,x_j)\le R\) for all \(i,j\in\{0,\dots,n\}\).  The
subcomplexes \(\abs{S(X)_R}\) are locally finite because bounded subsets
of~\(X\) are finite.  Up to quasi-isometry, there is a unique metric on
\(\abs{S(X)_R}\) for which the cells are uniformly bounded and the embedding
\(X\to\abs{S(X)_R}\) is a quasi-isometry.  This additional structure on
\(\abs{S(X)}\) is formalised in the following definition:

\begin{definition}
  \label{def:filtered_metric_CW}
  A \emph{filtered metric CW\brd{}complex} is a CW\brd{}complex~\(X\) together
  with a filtration \(X=\bigcup_{R\in\N} \Fil_R X\) and proper metrics~\(d_R\)
  on \(\Fil_R X\) that define the topology of \(\Fil_R X\) for all \(R\in\N\),
  such that each \(\Fil_RX\) is a locally finite, finite-dimensional
  subcomplex of~\(X\) with uniformly bounded cells and bounded boundary
  map~\(\delta\) and such that the embeddings \(\Fil_RX\to\Fil_{R+1}X\) are
  quasi\brd{}isometries for all \(R\in\N\).
  
  A morphism \(f\colon {X\to Y}\) between filtered metric CW\brd{}complexes is
  a cellular map with the property that for any \(R\in\N\) there exists
  \(R'\in\N\) such that \(f(\Fil_R X)\subseteq \Fil_{R'} Y\) and \(f|_{\Fil_R
    X}\colon \Fil_R X\to\Fil_{R'} Y\) is quasi\brd{}Lipschitz and bounded with
  respect to the norms defined in~\eqref{eq:cellular_norm}.
  
  If~\(X\) is a filtered metric CW\brd{}complex, then we turn \(X\times
  [0,1]\) into a filtered metric CW\brd{}complex using the standard cell
  decomposition, the filtration \((\Fil_R X)\times[0,1]\), and the metrics
  \(d_R\bigl((x,s),(y,t)\bigr)\defeq d_R(x,y)+\abs{t-s}\).  A \emph{homotopy}
  between two morphisms \(X\to Y\) is a morphism \(X\times[0,1]\to Y\).
\end{definition}

Up to isomorphism of filtered metric CW\brd{}complexes, only the equivalence
class of the filtration \((\Fil_R X)\) and the quasi\brd{}isometry classes of
the metrics~\(d_R\) matter.  Thus the filtered metric CW\brd{}complex
structure on \(\abs{S(X)}\) discussed above is unique up to isomorphism.

\begin{example}
  \label{exa:GCW_complex}
  Let~\(G\) be a finitely generated discrete group and let~\(X\) be a
  CW\brd{}complex equipped with a proper cellular action of~\(G\) and with
  only countably many cells.  Then~\(X\) becomes a filtered metric
  CW\brd{}complex in a canonical way.
  
  A \(G\)\nbd{}invariant subcomplex of~\(X\) is \(G\)\nbd{}compact if and only
  if the action of~\(G\) on its cells has only finitely many orbits.  Any
  finite subcomplex of~\(X\) is contained in a \(G\)\nbd{}invariant
  \(G\)\nbd{}compact subcomplex.  Since~\(X\) has countably many cells, there
  exists an increasing filtration \(X=\bigcup_{R\in\N} \Fil_RX\) where
  the~\(\Fil_RX\) are \(G\)\nbd{}invariant, \(G\)\nbd{}compact subcomplexes.
  This filtration is unique up to equivalence.  It follows from
  \(G\)\nbd{}compactness and properness that the subcomplexes \(\Fil_RX\) are
  locally finite and finite-dimensional and that the boundary map of
  \(\Fil_RX\) is bounded.
  
  Since the action of~\(G\) on \(\Fil_RX\) is proper and \(G\)\nbd{}compact,
  we can find a \(G\)\nbd{}invariant metric on \(\Fil_RX\) that defines the
  topology of \(\Fil_RX\), such that the maps \(G\to \Fil_RX\), \(g\mapsto
  gx\), are quasi\brd{}isometric for all \(x\in\Fil_RX\).  Any two such
  metrics are quasi\brd{}isometric.  Since the maps \(G\to \Fil_RX\),
  \(g\mapsto gx\), for some \(x\in\Fil_0X\) are quasi\brd{}isometric for all
  \(R\in\N\), it follows that the embeddings \(\Fil_RX\to\Fil_{R+1}X\) are
  also quasi\brd{}isometric.
  
  The above construction is natural, that is, \(G\)\nbd{}equivariant cellular
  maps become morphisms of filtered metric CW\brd{}complexes.  The same holds
  for \(G\)\nbd{}equivariant cellular homotopies.
\end{example}

Let~\(X\) be a filtered metric CW\brd{}complex and recall that~\(X_n\) denotes
its set of \(n\)\nbd{}cells.  It inherits a filtration \(X_n = \bigcup
\Fil_RX_n\).  Choose a point \(\xi_D\in D\) for each \(D\in X_n\) and view
this as a map \(X_n\to X\), \(D\mapsto \xi_D\).  This induces proper discrete
metrics on the sets \(\Fil_RX_n\) for all \(R,n\in\N\).  Our hypotheses
on~\(X\) ensure that the quasi\brd{}isometry classes of these metrics do not
depend on the choice of the points~\(\xi_D\).  Let
\[
\T\Cell_n(X)\defeq \varinjlim_{R\to\infty} \T(\Fil_RX_n),
\]
where~\(\T\) is \(\Hol\) or \(\Sch^k\) for \(k\in\R_+\cup\{\infty,\omega\}\).

\begin{lemma}
  \label{lem:tempered_chains_defined}
  Let~\(X\) be a filtered metric CW\brd{}complex.  Then the space
  \(\T\Cell_n(X)\) is independent of the auxiliary choices and the cellular
  boundary map extends to a bounded linear map \(\delta\colon
  \T\Cell_{n+1}(X)\to \T\Cell_n(X)\) for all \(n\in\N\).  Thus we obtain a
  bornological chain complex
  \(\T\Cell_\bullet(X)=(\T\Cell_n(X),\delta)_{n\ge0}\).  A morphism \(f\colon
  {X\to Y}\) induces a bounded chain map \(f_*\colon
  \T\Cell_\bullet(X)\to\T\Cell_\bullet(Y)\).

  A homotopy \(F\colon X\times[0,1]\to Y\) between two morphisms induces a
  bounded chain homotopy between the induced maps
  \(\T\Cell_\bullet(X)\to\T\Cell_\bullet(Y)\).
\end{lemma}

\begin{proof}
  The space \(\T(\Fil_RX_n)\) does not depend on the choice of the
  points~\(\xi_D\) because \(\T(\Fil_RX_n)\) only depends on the
  quasi\brd{}isometry class of the metric on~\(\Fil_RX_n\).  If \(\iota(\partial
  D,F)\neq0\) for some \(F\in\Fil_RX_n\), \(D\in\Fil_RX_{n+1}\), then \(F\cap
  \partial D\neq\emptyset\) and hence \(\abs{\ell(\xi_F)-\ell(\xi_D)}\)
  remains uniformly bounded.  Moreover, the number of summands in
  \(\delta(D)\) is uniformly bounded for \(D\in\Fil_RX_{n+1}\) because the
  boundary map~\(\delta\) on \(\Fil_RX\) is bounded.  If \(f\colon X\to Y\) is
  a morphism, then for each \(R\in\N\) there is \(R'\in\N\) such that
  \(f_*\colon \Cell_\bullet(X)\to\Cell_\bullet(Y)\) maps
  \(\Cell_\bullet(\Fil_RX)\to \Cell_\bullet(\Fil_{R'}Y)\).  If
  \(\iota(f,D,E)\neq0\) for two cells \(D\in\Fil_RX_n\), \(E\in\Fil_{R'}Y_n\),
  then \(f(D)\cap E\neq\emptyset\) and hence \(\ell(\xi_E) = O(\ell(\xi_D))\)
  because~\(f\) is quasi\brd{}Lipschitz and cells in \(\Fil_RX\) are uniformly
  bounded.  Since the sum \(\sum_E \abs{\iota(f,D,E)}\) is controlled
  uniformly for \(D\in\Fil_RX_n\), we obtain that~\(f_*\) extends to a bounded
  linear map \(f_*\colon \T\Cell_\bullet(X)\to\T\Cell_\bullet(Y)\).  The
  argument for homotopies is identical.
\end{proof}

Let~\(X\) be a discrete proper metric space and let \(S(X)\) be the associated
simplicial set.  Everything is set up so that \(\T\Cell_\bullet\abs{S(X)}\) is
equal to \(\T C_\bullet(X)\) if we turn \(\abs{S(X)}\) into a filtered metric
CW\brd{}complex as explained above.  If~\(Y\) is any filtered metric
CW\brd{}complex that is homotopy equivalent to \(\abs{S(X)}\) as such (it does
not suffice to have a homotopy equivalence just as topological spaces!), then
Lemma~\ref{lem:tempered_chains_defined} yields a bornological homotopy
equivalence between \(\T\Cell_\bullet(Y)\) and \(\T C_\bullet(X)\).  A
necessary condition for this is that~\(Y\) should be contractible and
\(\Fil_RY\) should be quasi\brd{}isometric to~\(X\), because \(\abs{S(X)}\)
certainly has these properties and they are invariant under homotopy
equivalence of filtered metric CW\brd{}complexes.  Thus we consider a filtered
metric CW\brd{}complex~\(Y\) together with a quasi\brd{}isometry \(X\to
\Fil_RY\) for some \(R\in\N\) in the following.

Recall that a CW\brd{}complex~\(Y\) is contractible if and only if all its
homotopy groups vanish.  We need a controlled version of this criterion:

\begin{definition}
  \label{def:uniformly_contractible}
  A filtered metric CW-complex~\(Y\) is \emph{uniformly contractible} if for
  any \(R\in\N\) there \(R'\in\N\) such that for any cellular map \(f\colon
  S^n\to \Fil_RY\) with \(n\le R\), \(\norm{f}_n\le R\), and diameter of
  \(f(S^n)\) at most~\(R\), there is an extension \(\bar{f}\colon
  D^{n+1}\to\Fil_{R'}Y\) with \(\norm{\bar{f}}_{n+1}\le R'\), and diameter of
  \(f(D^{n+1})\) at most~\(R'\).  Here we use the cell decomposition with
  minimal number of cells for the \(n\)\nbd{}sphere~\(S^n\) and the
  \({n+1}\)\brd{}cell~\(D^{n+1}\) (that is, \(D^{n+1}\) has~\(3\) cells, of
  dimension \(0\), \(n\) and \(n+1\), respectively, and the first two are
  those of \(S^n\subseteq D^{n+1}\)).
\end{definition}

\begin{theorem}
  \label{the:uniformly_contractible}
  Let~\(X\) be a discrete proper metric space and let~\(Y\) be a uniformly
  contractible, filtered metric CW\brd{}complex equipped with a
  quasi\brd{}isometry \(X\to\Fil_RY\) for some \(R\in\N\).  Then~\(Y\) is
  homotopy equivalent to \(\abs{S(X)}\) as a filtered metric CW\brd{}complex,
  and \(\T C_\bullet(X)\) is bornologically chain homotopy equivalent to
  \(\T\Cell_\bullet(Y)\).
\end{theorem}

\begin{proof}
  Let \(Z\) be another filtered metric CW\brd{}complex with a
  quasi\brd{}isometric map \(f\colon X\to\Fil_R Z\) for some \(R\in\N\).  We
  claim that there is a morphism of filtered metric CW\brd{}complexes
  \(f'\colon Z\to Y\) whose composition with~\(f\) is close to the given map
  \(X\to Y\); moreover, any two such morphisms are homotopic.  This universal
  property determines~\(Y\) uniquely up to homotopy equivalence of filtered
  metric CW\brd{}complexes.  It is easy to see that \(\abs{S(X)}\) is also
  uniformly contractible.  Hence the homotopy equivalence of \(\abs{S(X)}\)
  and~\(Y\) follows from the above claim.  The homotopy equivalence of
  \(\T\Cell_\bullet\abs{S(X)} = \T C_\bullet(X)\) then follows from
  Lemma~\ref{lem:tempered_chains_defined}.
  
  Thus it suffices to verify the universal property.  Both existence and
  uniqueness are proven by induction on skeleta.  The uniform contractibility
  guarantees the vanishing of the obstruction in the induction step.  We only
  write down the details for the existence part.  We have to construct
  functions \(\rho_n\colon \N\to\N\) and cellular maps \(f_n\colon Z^{(n)}\to
  Y\) with \(f_{n+1}|_{Z^{(n)}}=f_n\) for all \(n\in\N\), such that
  \(f_n(\Fil_RZ)\subseteq \Fil_{\rho_n(R)}(Y)\), the restriction of~\(f_n\) to
  \(\Fil_RZ^{(n)}\) is Lipschitz with Lipschitz constant \(\rho_n(R)\), and
  \(\norm{f_n|_{\Fil_RZ}}_n \le \rho_n(R)\), for all \(R\in\N\), \(n\in\N\).
  
  Since the given map \(X\to\Fil_RZ\subseteq \Fil_{R'}Z\) is a
  quasi\brd{}isometry for all \(R'\ge R\), we can find \(f_0\) and~\(\rho_0\)
  with the required properties.  In the induction step, we have to
  extend~\(f_n\) to \({n+1}\)\brd{}cells.  Such a cell~\(D\) is glued into
  \(Z^{(n)}\) by a characteristic map \(\chi_D\colon S^n\to Z^{(n)}\).  We
  have to extend \(f\circ\chi_D\colon S^n\to Y\) to a function \(D^{n+1}\to
  Y\).  If \(D\in\Fil_RZ_{n+1}\), then \(f\circ\chi_D\) takes values in
  \(\Fil_{\rho_n(R)}(Y)\).  It has Lipschitz constant
  \(C(R,n)\cdot\rho_n(R)\), where \(C(R,n)\) is the maximal diameter of
  \({n+1}\)\nbd{}cells in \(\Fil_RZ\), and it has \(n\)th norm at most
  \(\norm{\delta|_{\Fil_RZ}}_n\rho_n(R)\).  Therefore, uniform contractibility
  allows us to find \(\rho_{n+1}\colon \N\to\N\) and \(f_{n+1}\) with the
  required properties.
\end{proof}

\begin{example}
  \label{exa:uniformly_contractible_group}
  Let~\(G\) be a finitely generated discrete group and let~\(Y\) be a
  \(G\)\nbd{}CW-complex, equipped with the canonical filtered metric
  CW\brd{}complex structure defined in Example~\ref{exa:GCW_complex}.
  Then~\(Y\) is uniformly contractible if and only if it is contractible.  Of
  course, in this case we can get the assertion of
  Theorem~\ref{the:uniformly_contractible} much more easily from
  Proposition~\ref{pro:chains_resolutions}.
\end{example}

\begin{example}
  \label{exa:Rips}
  Let~\(X\) be a discrete proper metric space and let \(P_R(X)\), \(R\in\N\),
  be the associated Rips complexes.  The union \(P(X)=\bigcup P_R(X)\) is a
  filtered metric CW\brd{}complex if~\(X\) has bounded geometry.  Otherwise
  the complexes \(P_R(X)\) are not finite-dimensional and we have to modify
  the filtration slightly and take \(P(X)\defeq \bigcup P_R(X)^{(R)}\).  The
  finite dimension of \(\Fil_RY\) is needed in the proof of
  Theorem~\ref{the:uniformly_contractible}.  It is easy to check that \(P(X)\)
  is uniformly contractible.  Hence Theorem~\ref{the:uniformly_contractible}
  yields that \(\abs{S(X)}\) and \(P(X)\) are homotopy equivalent as filtered
  metric CW\brd{}complexes and that \(\T C_\bullet(X)\) and
  \(\T\Cell_\bullet(P(X))\) are bornologically homotopy equivalent.
  
  If the metric space~\(X\) is hyperbolic, then \(P_R(X)\) is uniformly
  contractible for sufficiently large~\(R\).  Thus \(\abs{S(X)}\) is homotopy
  equivalent to \(P_R(X)\) as a filtered metric CW\brd{}complex and \(\T
  C_\bullet(X)\cong \T\Cell_\bullet(P_R(X))\).  Since hyperbolic spaces have
  combings of linear growth, Theorem~\ref{the:combable_contract} implies that
  the complexes \(\Sch\Cell_\bullet(P_R(X))\),
  \(\Sch^\omega\Cell_\bullet(P_R(X))\), and \(\Hol\Cell_\bullet(P_R(X))\) are
  bornologically contractible.
\end{example}


\bib{Allcock-Gersten}{article}{
  author={Allcock, Daniel J.},
  author={Gersten, Stephen M.},
  title={A homological characterization of hyperbolic groups},
  journal={Invent. Math.},
  volume={135},
  date={1999},
  number={3},
  pages={723--742},
  issn={0020-9910},
  review={\MR {1669272}},
  doi={10.1007/s002220050299},
}

\bib{Bridson:Combings_semidirect}{article}{
  author={Bridson, Martin R.},
  title={Combings of semidirect products and $3$\nobreakdash -manifold groups},
  journal={Geom. Funct. Anal.},
  volume={3},
  date={1993},
  number={3},
  pages={263--278},
  issn={1016-443X},
  review={\MR {1215781}},
  doi={10.1007/BF01895689},
}

\bib{Connes-Moscovici:Novikov_Hyperbolic}{article}{
  author={Connes, Alain},
  author={Moscovici, Henri},
  title={Cyclic cohomology, the Novikov conjecture and hyperbolic groups},
  journal={Topology},
  volume={29},
  date={1990},
  number={3},
  pages={345--388},
  issn={0040-9383},
  review={\MR {1066176}},
  doi={10.1016/0040-9383(90)90003-3},
}

\bib{Epstein:Automatic}{book}{
  author={Cannon, James W.},
  author={Epstein, David B. A.},
  author={Holt, Derek F.},
  author={Levy, Silvio V. F.},
  author={Paterson, Michael S.},
  author={Thurston, William P.},
  title={Word processing in groups},
  publisher={Jones and Bartlett Publishers},
  place={Boston, MA},
  date={1992},
  pages={xii+330},
  isbn={0-86720-244-0},
  review={\MR {1161694}},
}

\bib{Gilman-Holt-Rees}{article}{
  author={Gilman, Robert H.},
  author={Holt, Derek F.},
  author={Rees, Sarah},
  title={Combing nilpotent and polycyclic groups},
  journal={Internat. J. Algebra Comput.},
  volume={9},
  date={1999},
  number={2},
  pages={135--155},
  issn={0218-1967},
  review={\MR {1703070}},
  doi={10.1142/S0218196799000102},
}

\bib{Hogbe-Nlend:Completions}{article}{
  author={Hogbe-Nlend, Henri},
  title={Compl\'etion, tenseurs et nucl\'earit\'e en bornologie},
  journal={J. Math. Pures Appl. (9)},
  volume={49},
  date={1970},
  pages={193--288},
  review={\MR {0279557}},
}

\bib{Hogbe-Nlend:Bornologies}{book}{
  author={Hogbe-Nlend, Henri},
  title={Bornologies and functional analysis},
  publisher={North-Holland Publishing Co.},
  place={Amsterdam},
  date={1977},
  pages={xii+144},
  isbn={0-7204-0712-5},
  review={\MR {0500064}},
}

\bib{Ji:Schwartz_cohomology}{article}{
  author={Ji, Ronghui},
  title={Smooth dense subalgebras of reduced group $C^*$-algebras, Schwartz cohomology of groups, and cyclic cohomology},
  journal={J. Funct. Anal.},
  volume={107},
  date={1992},
  number={1},
  pages={1--33},
  issn={0022-1236},
  review={\MR {1165864}},
  doi={10.1016/0022-1236(92)90098-4},
}

\bib{Meyer:Born_Top}{article}{
  author={Meyer, Ralf},
  title={Bornological versus topological analysis in metrizable spaces},
  conference={ title={Banach algebras and their applications}, place={Edmonton}, },
  book={ series={Contemp. Math.}, volume={363}, publisher={Amer. Math. Soc.}, place={Providence, RI}, },
  date={2004},
  pages={249--278},
  review={\MR {2097966}},
  doi={10.1090/conm/363},
}

\bib{Meyer:Embed_derived}{article}{
  author={Meyer, Ralf},
  title={Embeddings of derived categories of bornological modules},
  date={2004},
  status={eprint},
  note={\arxiv {math.FA/0410596}},
}

\bib{Meyer:Poly_Oberwolfach}{article}{
  author={Meyer, Ralf},
  title={Polynomial growth cohomology for combable groups},
  volume={45},
  journal={Oberwolfach reports},
  issn={1660-8933},
  date={2004},
  pages={2380--2384},
}

\bib{Meyer:Homological_Schwartz}{article}{
  author={Meyer, Ralf},
  title={Homological algebra for Schwartz algebras of reductive \(p\)\nobreakdash -adic groups},
  conference={ title={Noncommutative geometry and number theory}, },
  book={ series={Aspects Math., E37}, publisher={Vieweg}, place={Wiesbaden}, },
  date={2006},
  pages={263--300},
  review={\MR {2327309}},
  doi={10.1007/978-3-8348-0352-8_13},
}

\bib{Mineyev:Homology_Combable}{article}{
  author={Mineyev, Igor},
  title={$l_1$-homology of combable groups and $3$-manifold groups},
  note={International Conference on Geometric and Combinatorial Methods in Group Theory and Semigroup Theory (Lincoln, NE, 2000)},
  journal={Internat. J. Algebra Comput.},
  volume={12},
  date={2002},
  number={1-2},
  pages={341--355},
  issn={0218-1967},
  review={\MR {1902370}},
  doi={10.1142/S0218196702000869},
}



\begin{bibdiv}
  \begin{biblist}
    \bibselect{references}
  \end{biblist}
\end{bibdiv}
\end{document}